\documentclass{siamltex}
\usepackage{amsfonts}        %
\usepackage{ipe}
\usepackage{epsfig}
\usepackage{amsmath,amssymb} %
\usepackage{latexsym}        %
\usepackage{t1enc}

\def\Bbb{\mathbb}

\def\RR#1{{\Bbb R}^{#1}}
\def\RNUM{{\Bbb R}}

\def\CNUM{{\Bbb C}}

\def\NORM#1#2{|\!| {#1} {|\!|}_{{#2}}}

\def\d{\,d}

\def\EXP#1{e^{#1}}
\def\DET{\det}

\def\GRAD{\nabla}
\def\LAP{\Delta}

\def\PD#1#2{\frac{\partial{#1}}{\partial{#2}}}
\def\DF#1#2{\frac{{\mathrm d}{#1}}{{\mathrm d}{#2}}}
\def\DDF#1#2{\frac{{\mathrm d}^2{#1}}{{\mathrm d}{#2}^2}}

\def\PDER#1#2{{\mathrm{D}_{#1}{#2}}}
\def\DFUN#1{{{#1}'}}

\def\SPACE{\;\;\;}

\def\COMMA{\,,}
\def\PERIOD{\,.}
\def\SEP{{\,|\,}}
\def\VIZ#1{(\ref{#1})}

\def\num#1{{\tt{#1}}}

\def\Im{i}
\def\REAL{\mathrm{Re}\,}
\def\IMAG{\mathrm{Im}\,}
\def\ARG{\mathrm{arg}\,}
\def\BIGO{\mathrm{O}}
\def\CONST{\mathrm{const}}
\def\qed{$\Box$}
\def\PROOF{{\sc Proof:}\ }

\def\ASYM#1{{#1}^{-1/\sigma-\Im\omega/\kappa}}
\def\BVPQ{\VIZ{clg3}-\VIZ{clg4}}
\def\BVPS{\VIZ{clg18}-\VIZ{clg18b}}
\def\NU{\kappa, \omega}

\def\ALPHA{\epsilon,\delta}
\def\JFUN{\mathcal{J}}
\def\LOPER{\mathcal{L}}

\def\BSPT{\mathcal{X}_{\vartheta}}
\def\BSPB{\mathcal{B}}
\def\PER{[0,2\pi)}
\def\WSP{\mathcal{W}}
\def\ZSP{\mathcal{Z}}
\def\LAM{\Lambda}
\def\LAMT{\tilde\Lambda}
\def\lambdat{{\tilde\lambda}}
\def\ZFUN{Z}
\def\Mm{\mathcal{M}}
\def\Ll{\mathcal{L}}

\def\Pp{\mathcal{P}}
\def\Ss{\mathcal{S}}

\newtheorem{remark}{Remark}[section]

\title{%
On self-similar singular solutions of the complex Ginzburg-Landau
equation.
}%
\author{%
Petr Plech\'a\v{c}\thanks{Department of Mathematical Sciences, University of Delaware,
  Newark, DE 19716, e-mail: {\tt plechac@math.udel.edu}.  The research of
  PP was supported in part by grant UDRF-19990415 from University of Delaware
  Research Foundation.}
\and
Vladim\'\i{}r \v{S}ver\'ak\thanks{School of Mathematics, 
University of Minnesota, Minneapolis, MN 55455, e-mail: 
{\tt sverak@math.umn.edu}.
The research of VS was supported in part by grant  DMS-9877055 from the National Science
Foundation.}
}
\begin{document}

\maketitle

\begin{abstract}
We address the open problem of  existence of singularities for
the complex Ginzburg-Landau equation. Using a combination of rigourous results 
and numerical computations, we  describe a countable family 
of self-similar singularities. Our analysis includes the  super-critical non-linear 
Schr\"odinger equation as a special case,
and most of the described singularities are new even
in that situation.
We also consider  the problem of stability of these singularities.
\end{abstract}

\section{Introduction}

In this paper we study singular solutions to
the initial value problem for the complex Ginzburg-Landau 
equation (CGL)
\begin{eqnarray}\label{clg1}
\Im \PD{u}{t} + (1-\Im\epsilon)\LAP u + 
   (1+\Im\delta) |u|^{2\sigma} u & = & f\COMMA
   \SPACE\mbox{in $\RR{d}\times (0,T)$}\COMMA \label{clg1a} \\
  u(x,0) & = & u_0(x)\PERIOD \label{clg1b}
\end{eqnarray}
where $u=u(x,t)$ is a complex-valued function defined 
in $\RR{d}\times (0,T)$ and satisfying suitable decay conditions as
$|x|\to\infty$,
 the parameters $\epsilon$, $\delta$, $\sigma$
are non-negative real numbers,  $u_0$ is a given initial condition,
which is assumed
to be smooth, with suitable decay as $|x|\to \infty$, and $f=f(x,t)$ is a
given function, also assumed to be smooth, with suitable decay 
as $|x|\to \infty$.

We are mostly interested in the case $\epsilon>0$, $\delta\geq 0$, 
and $2/d <\sigma <2/d+1/2$.
However, some of our results are new even
for the case $\epsilon\geq 0$, $\delta\geq 0$, $\sigma>2/d$, i.e.,
they also include the super-critical non-linear Schr\"odinger equation.
Under the assumption $\epsilon>0$, $0<\sigma <2/d+1/2$
it is possible to prove that for each smooth $u_0$ and $f$, with an
appropriate decay at infinity, the problem \VIZ{clg1a}-\VIZ{clg1b}
 has a {\it suitable weak solution}, see \cite{DoeringGibbonLevermore} and \cite{XYan}.
Such a solution is regular away from a closed set  
$\Ss\subset \RR{d}\times (0,T)$ with $\Pp^{d-2/\sigma}(\Ss)=0$,
where $\Pp^\alpha$ denotes the {\it parabolic $\alpha$-dimensional
Hausdorff measure}, see \cite{XYan}.
It has been an open problem whether the singular set $\Ss$ can be
non-empty.

We present a very strong evidence, which is
based on a combination of  rigorous analysis and  numerical computations,
that 
singularities may indeed exist. We shall see that there appears to exist a
countable set of different types of singular solutions. Among these
solutions
we identify (numerically) those which are stable under a suitable notion of stability 
defined in Section~\ref{sec4}. It turns out that while most of the
solutions are unstable, in certain cases there may exist more than one type of
stable singularities.

The case of the non-linear Schr\"odinger equation (NLS)
(i.e., $\epsilon=\delta=0$) and $f=0$ has been studied 
by  numerous authors, see for example
\cite{McLaughlinPapanicolaou,KopellLandman} or the monograph 
\cite{SulemSulem}. 
In particular, Zacharov (\cite{Zacharov1}) conjectured the existence of 
self-similar singularities of the form
\begin{equation} \label{selfsim}
u(x,t) = \left( 2\kappa
(T-t)\right)^{-\frac{1}{2}(\frac{1}{\sigma}+\Im\frac{\omega}{\kappa})}
 \,Q\left((2\kappa(T-t))^{-1/2} |x| \right)\COMMA
\end{equation}
where $Q(\xi)$ is a complex valued function defined on 
 $(0,\infty)$, with asymptotic behavior
$$
Q(\xi) \sim \ASYM{\xi}\COMMA\SPACE
\mbox{as} \SPACE \xi\to\infty \PERIOD
$$

While no rigorous proof of this conjecture seems to be
available, there is an overwhelming evidence based on numerical and
formal analytical calculations supporting the existence of such
singularities
(see, e.g., \cite{LeMesurier,LeMesurierPapanicolau, SulemSulem}).
We also refer the reader to \cite{KopellLandman} for a rigorous result
supporting the conjecture.

\medskip

In this paper we will argue that these singularities persist also
for a certain range of $\epsilon>0$ and $\delta>0$. In fact, we shall
find many new self-similar singularities even for the case $\epsilon=
\delta=0$ (and $f=0$). Using the self-similar singular solutions of the
form \VIZ{selfsim} one can easily construct singular solutions of
\VIZ{clg1a} and \VIZ{clg1b} with compactly supported, smooth $u_0$ and 
compactly supported, smooth $f$.

\medskip

From the form \VIZ{selfsim} of self-similar solutions one obtains the 
following boundary value problem for the function $Q$:
\begin{eqnarray}
&& (1-\Im\epsilon)(Q'' + \frac{d-1}{\xi}Q') +
      \Im\kappa\xi Q' + 
      \Im\frac{\kappa}{\sigma}Q - \omega Q + 
      (1+\Im\delta)|Q|^{2\sigma}Q = 0 \COMMA \label{clg3} \\
&& Q'(0)=0 \COMMA\;\;\mbox{and} \label{clg4a}\\
&& Q(\xi)\sim \ASYM{\xi}\;\;\mbox{as}\;\; \xi\to\infty\PERIOD \label{clg4}
\end{eqnarray}

In general, this problem does not have a solution for arbitrary values
of parameters $\kappa$, $\omega$ and hence the unknowns in \BVPQ\  
are $Q$, $\kappa$ and $\omega$.

When $\epsilon=0$, see, for example, \cite{SulemSulem} for a discussion of
\VIZ{clg4}.
For $\epsilon >0$ the condition $Q(\xi)\sim\ASYM{\xi}$ at  infinity is
dictated by the
partial regularity result mentioned earlier, since $u$ defined by
\VIZ{selfsim} must
be regular at almost all points of the form $(x,T)$, $x\ne 0$.
The presented results are based on a detailed analysis of the boundary 
value problem \BVPQ.  The analysis employs both analytical and
numerical techniques and is naturally divided into the following  
two steps:

\begin{description}
\item[{\rm (i)}] First, we rigorously prove that for $\xi_1\ge 1$ and
   for each set of
  parameter values $\epsilon\geq 0$, $\delta\geq 0$ and $2/d < \sigma
  < 2/d + 1/2$ there exists a two-dimensional manifold 
  of solutions to \VIZ{clg3} on $(\xi_1,\infty)$ with the correct
  asymptotic behavior at infinity. 
\item[\rm(ii)] Second, we solve numerically the boundary value problem 
  on $(0,\xi_1)$ with an appropriate (approximate) boundary condition for 
  $Q$ at the  boundary point $\xi=\xi_1$.  
  In this step the choice of $\xi_1$ is based on numerical evidence 
  of convergence, 
  and not on the rigorous estimates obtained in (i),
  which contain constants we did not try to evaluate exactly.
\end{description}

\medskip

We briefly summarize main 
results of our calculations.  We note a (scaling) symmetry in the
problem \BVPQ: If $\lambda>0$, $\theta\in\RNUM$ and
$(Q(\xi),\kappa,\omega)$ is a solution of \BVPQ\
then
 $(\lambda^{1/\sigma + \Im\theta} Q(\lambda\xi), 
\lambda^2\kappa,\lambda^2\omega)$ is also a solution.
We chose a representative in each of these families of solutions by
imposing suitable normalization conditions. We mostly work with the
normalization $\omega=1$ and $\IMAG{Q(0)}=0$, $Q(0)>0$. Such solutions
will be
called {\it normalized} solutions and they are uniquely determined by two
parameters $(\kappa,\mu)$, where $\mu=Q(0)$. Some numerical results
are better presented in a different parameterization in which $Q(0)=1$ is
fixed,
$(\kappa,\omega)$ are used as parameters. We shall state the use of the
latter
normalization explicitly whenever it is used. 

We observe that for fixed values of $d\geq 1$, $\sigma>2/d$ and
$\epsilon=\delta=0$ there exists a countable family of normalized
solutions $(Q_j(\xi),\kappa_j)=(Q_j(\xi),\kappa_j,1)$, $j=1,2,\dots$ of \BVPQ. 
The $j-$th profile $|Q_j|$, when extended to $(-\infty,\infty)$ as an even function, 
has exactly
$j$ local maxima, and somewhat resembles the profile of the $j-$th state 
of an elementary quantum mechanical oscillator. A precise quantum mechanical interpretation
of the solutions is more complicated, and is related to the so-called {\it resonances}
or {\it quasi-stationary states}.
The first solution of the family has 
been known, see \cite{McLaughlinPapanicolaou} or
\cite{SulemSulem}, for example. We are not aware of any mentioning of the other
solutions in the literature.
All the solutions persist if $\epsilon$ and $\delta$ are
perturbed to (small) strictly positive values. 
We now describe the behavior of these perturbed
solutions for $\delta=0$. We let $\mu_j=Q_j(0)$. As mentioned above,
the solution $(Q_j,\kappa_j)$ is determined by $(\kappa_j,\mu_j)$.
 We observe that a branch of solutions
parameterized by $(\kappa_j(\epsilon),\mu_j(\epsilon))$ emanates
from each point $(\kappa_j,\mu_j)$.  We plot the curves $(\epsilon,\kappa_j(\epsilon))$,
$j=1,\dots, 5$,
in Figure~\ref{fig1} for $d=1$, $\sigma=2.3$ and in Figure~\ref{fig2}
for $d=3$, $\sigma=1$.
In the other figures we  plot the profiles $|Q^\epsilon_j(\xi)|$ of the corresponding
solutions $Q^\epsilon_j(\xi)$ at certain points along
each branch. 

An interesting feature observed in the behavior of the
branches, is the existence of a turning point on each branch at
$\epsilon=\epsilon^*_j$. We see from the graphs that as $\epsilon$ 
returns to zero along 
the branch, $\kappa_j(\epsilon)$ tends to zero, suggesting that the
solution $Q^\epsilon_j$ converges to a (radial) solution of the
equation
\begin{equation}\label{clg5}
\LAP Q - Q + |Q|^{2\sigma}Q = 0\COMMA\;\;\mbox{in $\RR{d}$}\COMMA
\end{equation}
satisfying $Q(x)\to 0$ as $|x|\to \infty$. 

Our computations presented in Section~\ref{sec3} 
clarify the structure of the diagram which turns out to be slightly more
complicated than the picture suggested above. We conjecture the following:
If $d>1$, the solutions $Q^\epsilon_j$  corresponding to branch 
$(\kappa_j(\epsilon),\mu_j(\epsilon))$ with an odd
index $j=2k-1$  converge (as $\epsilon$ and $\kappa_j(\epsilon)$ approach zero)
to the $k$-th (normalized) radial solution of \VIZ{clg5}. We recall that the first
of these solutions is usually called {\it the  ground state}, and that for $d=1$ there
are no other solutions of \VIZ{clg5} satisfying the appropriate boundary conditions.
(See, for example, \cite{SulemSulem} for more details.) 
If $d=1$ and also
for $j=2k$ in the case $d>1$, as $\epsilon$ and $\kappa_j(\epsilon)$ approach zero,
 the profiles $Q^\epsilon_j$ separate into $j$ approximate ground
states which move away from each other. In particular, for $j=2k$ the profiles
converge locally uniformly to zero. When $d=1$ and $j=2k-1$, the profiles
converge locally uniformly to the ground state. 
  
We tabulate results of our numerical calculations in Table~\ref{tab1} and
Table~\ref{tab2}.

We also looked at branches of solutions when  $\delta$ is related to
$\epsilon$ by $\delta=r\epsilon$, with $r>0$ of order $10^{-1}$
and $10^0$. The behavior was similar to the case $\delta=0$, with
the turning point $\epsilon^*$ getting closer to zero as $r$ increased,
as one might heuristically expect.

\medskip

Questions related to stability of the singularities are addressed in Section~\ref{sec4}.
Our calculations indicate that for the non-linear Schr\"odinger equation all the new singularities
we found are unstable, and the singularity corresponding to $(\kappa_1,\mu_1)$ 
is stable. The situation is more complicated for $\epsilon>0$, see Section~\ref{sec4}
for details.

\medskip

Our interest in singular solutions to CGL
stems from analogies between \VIZ{clg1} in the case $d=3$,
$\sigma=1$, $\epsilon,\delta>0$ and the three-dimensional 
Navier-Stokes equation (NSE).
The two equations have the same scaling properties and the same energy
identity.
Moreover, the existence and partial regularity theory of weak solutions
for
NSE and CGL are similar (with CGL being technically easier), see
\cite{XYan}.
The analogy between NSE and CGL may be rather
superficial and may break down at any deeper level. However, at the same
time there are
no known properties of solutions to the Navier-Stokes equation which would
prevent the
same scenario as presented here for CGL.

The formula for the Navier-Stokes equation corresponding to \VIZ{selfsim}
would be
\begin{equation}\label{selfsimnse}
u(x,t) \sim (2\kappa(T-t))^{-1/2} 
        \, U\left((2\kappa(T-t))^{-1/2} x,\tau \right)\COMMA
\end{equation}
where 
$\tau=\frac{1}{2\kappa}\ln\frac{T}{T-t}$,
and $U$ is a suitable divergence-free vector field periodic in $\tau$ with
suitable
decay in the self-similar variable $y=(2\kappa(T-t))^{-1/2}x$.

The case, $\partial U/\partial\tau \equiv 0$, 
was already considered by Leray \cite{Leray}. It was proved in 
\cite{NecasRuzickaSverak}
and in greater generality in \cite{Tsai} that NSE does not admit
non-trivial solutions
of the form \VIZ{selfsimnse} with $U$ independent of $\tau$.
The problem is open for  $U$ periodic in $\tau$.

\medskip 

We finish the introduction with the following speculation. Most of the 
singularities we have found are unstable, hence it is unlikely they would
be observed in direct numerical simulations of the initial value
problem \VIZ{clg1}-\VIZ{clg1b} or in physical experiments that are modeled
by
CGL. Could it perhaps be
the case that NSE does admit singular solutions (say of the form
\VIZ{selfsimnse}), but {\it all} of them are unstable and therefore more
or
less impossible to be detected in direct numerical simulations or physical
experiments ?

This intriguing scenario was once suggested to one of the authors by
\mbox{Sergiu}  \mbox{Klainerman} during a lunch-break conversation at a conference in 
Southern California.

\section{Analysis of the profile equation at infinity}\label{sec2}

In this section we study solutions of \VIZ{clg3} in the interval
$(\xi_1,\infty)$ satisfying the condition $Q(\xi)\sim \ASYM{\xi}$
as $\xi\to\infty$. Heuristically one expects that the behavior of
such solutions is mainly governed by the linear part of
\VIZ{clg3}:
\begin{equation}\label{clgL}
(1-\Im\epsilon) u'' + (1-\Im\epsilon) \frac{d-1}{\xi} u' +
\Im\kappa\xi u' + \frac{\Im\kappa}{\sigma} u - \omega u = 0
\end{equation}
Equation \VIZ{clgL} is equivalent to  Kummer's equation, also
known as the confluent hypergeometric equation. The solutions of this
equation are well-understood, see for example \cite{Slater}, and one can 
hence get 
a more or less complete picture of the behavior of solutions to
\VIZ{clgL}.  We shall use analytical tools from the theory of
confluent hyper-geometric equations to describe solutions of the full
equation \VIZ{clg3}.

A canonical form of Kummer's equation is
\begin{equation}\label{clgK}
z\DDF{w}{z} + (b-z) \DF{w}{z} - a w = 0\COMMA
\end{equation}
and the equation \VIZ{clgL} is transformed into this form by letting
\begin{equation}\label{clg6}
z=\frac{-\Im\kappa}{(1-\Im\epsilon)}\frac{\xi^2}{2}\COMMA\SPACE\SPACE
a=\frac{1}{2}\left(\frac{1}{\sigma}+\frac{\Im\omega}{\kappa}\right)\COMMA
\SPACE\SPACE
b=\frac{d}{2}\PERIOD
\end{equation}
There is voluminous literature on this equation and properties of special functions
(confluent hyper-geometric functions) which appear as its solutions. 
We recall some of the properties and, for the convenience of the
reader, we also sketch how to derive them.
For more details about confluent hyper-geometric
functions we refer the reader to
\cite{Slater}, \cite{Luke}, \cite{Whittaker}.

\medskip

A classical formula for a solution of \VIZ{clgK} is given by 
\begin{equation}\label{clg7}
U(a,b,z) = \frac{1}{\Gamma(a)}\int_0^\infty \EXP{-tz} t^{a-1}
(1+t)^{b-a-1}\,dt\COMMA
\end{equation}
where $\Gamma$ is Euler's gamma function. The integral is clearly
well-defined for $\REAL{a}>0$ and $\REAL{z}>0$. The factor
$1/\Gamma(a)$ is not essential for our analysis, nevertheless,
we include it to keep our notation in agreement with the standard
one. The role of this factor is to normalize the leading term in the asymptotic
$U(a,b,z)=z^{-a}(1+\BIGO(z^{-1}))$ as $z\to\infty$.
It is easy to check by direct calculation that the function $U$ given by \VIZ{clg7} solves
the equation \VIZ{clgK}. We have
$$
\frac{d^k}{dz^k} U(a,b,z) = \frac{1}{\Gamma(a)}\int_0^\infty (-t)^k \EXP{-tz} t^{a-1}
(1+t)^{b-a-1}\,dt\COMMA 
$$
and after substituting to \VIZ{clgK} we see from a 
simple integration by parts that the equation is satisfied.

By letting $zt=s$ in \VIZ{clg7} we obtain (for $\REAL{z}>0$,
$\REAL{a}>0$)
\begin{equation}\label{clg8}
U(a,b,z) = \frac{1}{\Gamma(a)} z^{-a}\int_0^\infty \EXP{-s} s^{a-1}
\left(1+\frac{s}{z}\right)^{b-a-1}\,ds\PERIOD
\end{equation}
The above expression is used to extend the definition of $U$ for
$\REAL{a}>0$ and $z\in\CNUM$ with $-\pi < \ARG{z} < \pi$, since the
integral is convergent and an analytic function of $z$ under these assumptions.
By formally expanding the term $(1+s/z)^{b-a-1}$ and integrating the
resulting (formal) series term-by-term, we obtain the following
asymptotic expansion
\begin{equation}\label{clg9}
U(a,b,z) = z^{-a} \left( \sum_{k=0}^{n-1} \frac{(a)_k(1+a-b)_k}{k!} (-z)^{-k}
                          + \BIGO(|z|^{-n}) \right) \COMMA
\end{equation}
where $(a)_0=1$, $(a)_k=a(a+1)\dots(a+k-1)$. The formal
calculations of the asymptotic expansion \VIZ{clg9} can be easily justified
rigorously by splitting the
integral in \VIZ{clg8} as $\int_0^\infty = \int_0^{s_1} +
\int_{s_1}^\infty$ with $s_1 = |z|^{1/2}$, for example.

One can  check easily by a direct calculation that the function 
\begin{equation}\label{clgV}
V(a,b,z) = \EXP{z} U(b-a,b,-z)
\end{equation}
is another solution of Kummer's equation and that the functions
$U$ and $V$ are linearly independent. The Wronskian of $U$ and $V$
is given by
\begin{equation}\label{clg10}
U\DF{V}{z} - V \DF{U}{z} = \EXP{\pm\Im\pi (b-a)} z^{-b} \EXP{z}\COMMA
\end{equation}
where the sign $+$ is for $\IMAG{z}>0$ and $-$ in the opposite case.
(Formula \VIZ{clg10} is easily derived from the fact that the Wronskian
satisfies the differential equation $y'+(b/z-1)y=0$ and from
the asymptotic expansion \VIZ{clg9}.)

We need to know the behavior of $U$ and $V$ in the region $-\pi/2\leq
\ARG{z}<0$. Taking into account the definition of $V$, we see
that it is sufficient to control $U$ in the region
$-\pi/2\leq\ARG{z}<\pi$. Formula \VIZ{clg8} is suitable for
analysis in this region if $\REAL{(b-a)>0}$ since the integral in
\VIZ{clg8} is the uniformly absolutely convergent whenever $z$
approaches a point in $(-\infty,0)$ from the upper half-plane. In our
applications the condition $\REAL{(b-a)>0}$ is always satisfied and 
therefore \VIZ{clg8} is sufficient for our analysis.

We now have sufficient information about the solutions of \VIZ{clgK},
and hence also \VIZ{clgL}, to be able to proceed with the analysis of the
inhomogeneous equation
\begin{equation}\label{clgLI}
(1-\Im\epsilon) u'' + (1-\Im\epsilon) \frac{d-1}{\xi} u' +
\Im\kappa\xi u' + \frac{\Im\kappa}{\sigma} u - \omega u = f(\xi)\COMMA
\end{equation}
for $\xi\in(\xi_1,\infty)$. We assume that the function $f$ is decaying sufficiently
fast as $\xi\to\infty$. We are interested in solutions of \VIZ{clgLI}
which have the asymptotics $u\sim \ASYM{\xi}$ as
$\xi\to\infty$. We denote $P$, $E$ two linearly independent solutions of \VIZ{clg3}
\begin{eqnarray*}
&& P(\xi) \equiv P(\kappa,\omega,\epsilon;\xi) = U(a,b,z)\COMMA \\
&& E(\xi) \equiv E(\kappa,\omega,\epsilon;\xi) = V(a,b,z)\COMMA \\
&&\mbox{where }\;
z=\frac{-\Im\kappa}{1-\Im\epsilon}\frac{\xi^2}{2}\COMMA\; a=
\frac{1}{2}\left(\frac{1}{\sigma}+\Im\frac{\omega}{\kappa}\right)\COMMA\; 
b=\frac{d}{2}\PERIOD
\end{eqnarray*}
The parameters $d$ and $\sigma$ are held fixed in the perturbation
analysis described below, therefore we do not indicate the dependence
of $P$ and $E$ on them. The Wronskian $W=PE'-P'E$ is easily
computed from 
\VIZ{clg10}.
\begin{equation}\label{clg11}
W(\xi)= W(\omega,\kappa,\epsilon;\xi) =
\frac{-\Im\kappa}{1-\Im\epsilon}  \EXP{\pm\pi\Im(b-a)} \xi z^{-b} \EXP{z}\PERIOD
\end{equation}

The next step is to use the standard variation of constant to obtain
solutions to \VIZ{clgLI} satisfying $u(\xi)\sim\ASYM{\xi}$ as
$\xi\to\infty$. 
We write the solution in the form
\begin{eqnarray}\label{clg12}
u(\xi) & = & c_1(\xi) P(\xi) + c_2(\xi) E(\xi)\COMMA\;\mbox{with} \\
u'(\xi) & = & c_1(\xi) P'(\xi) + c_2(\xi) E'(\xi)\COMMA
\end{eqnarray}
and $c_2(\xi)\to 0$ ``sufficiently fast'' as $\xi\to\infty$. We
obtain 
\begin{equation}\label{clg13}
c_1'(\xi) = -f(\xi)\frac{E(\xi)}{(1-\Im\epsilon)W(\xi)}\COMMA\;\;\;
c_2'(\xi) =  f(\xi)\frac{P(\xi)}{(1-\Im\epsilon)W(\xi)}\COMMA
\end{equation}
which together with the condition $c_2(\xi)\to 0$ at the infinity
gives a formal expression for the solution
\begin{equation}\label{clg14}
u(\xi) = \gamma P(\xi) + \int_{\xi_1}^\infty K(\xi,\eta) f(\eta)
\d\eta\COMMA
\end{equation}
where $\gamma\in\CNUM$ is a constant and
$$
K(\xi,\eta) = \left\{ \begin{array}{ll}
              -\frac{1}{(1-\Im\epsilon)}P(\xi)E(\eta)W^{-1}(\eta) &
            \mbox{for } \xi_1 < \eta \leq \xi \\
              -\frac{1}{(1-\Im\epsilon)}E(\xi)P(\eta)W^{-1}(\eta) &
            \mbox{for } \xi  \leq \eta
                      \end{array} \right.\PERIOD
$$
The strategy for finding solutions
of the non-linear problem \VIZ{clg3} in $(\xi_1,\infty)$ 
with the required decay $u(\xi)\sim
\ASYM{\xi}$ at infinity should now be clear. Heuristically we expect that the
manifold of such solutions will be a  deformation of the
one-dimensional complex subspace $\{\gamma P\SEP \gamma\in\CNUM\}$, at 
least in a neighborhood of the origin. The deformation is from $\gamma P$ to 
the fixed point of the operator $T$
$$
u \mapsto T u(\xi) = \gamma P(\xi) -\int_{\xi_1}^\infty (1+i\delta)
                        K(\xi,\eta)\,
                    |u(\eta)|^{2\sigma}u(\eta)\d\eta \PERIOD
$$

The outlined strategy can be successfully carried out by using the
properties of Kummer's functions recalled above. The fixed point
theorem can be applied in the Banach space 
$$
\BSPT=\{ u\in C([\xi_1,\infty))\SEP \sup_{\xi\geq\xi_1}
|\xi|^{1/\sigma-\vartheta}|u(\xi)| <\infty\}
$$
equipped with the norm
$$
\NORM{u}{\vartheta} = \sup_{\xi\geq\xi_1} |\xi|^{1/\sigma-\vartheta}|u(\xi)|\PERIOD
$$

This approach is obviously standard. However, there are
some subtle points in the situation studied here due to oscillatory
behavior of the function $E$, see the Appendix.

The main result of this section, which can be 
derived in a fully rigorous way from the above analysis is the following
theorem.
(A complete proof of the theorem is presented in the Appendix.)

\begin{theorem}\label{thm1}
Assume $1\leq d \leq 3$, $2/d<\sigma < 2/d+1/2$, 
$0<\kappa_1<\kappa_2$, $\omega_1<\omega_2$, $0<\epsilon_1$, 
$\delta_1<\delta_2$.

There exists $\rho_1>0$ such that for each $\xi_1\ge 1$ and each
$$
(\beta,\kappa,\omega,\epsilon,\delta)\in \LAM_{\rho_1}\equiv 
\{\beta\in\CNUM\SEP|\beta|\leq\rho_1\}\times
[\kappa_1,\kappa_2] \times
[\omega_1,\omega_2] \times
[0,\epsilon_1] \times
[\delta_1,\delta_2]
$$
the boundary value problem
\begin{eqnarray}
&& (1-\Im\epsilon)(Q'' +  \frac{d-1}{\xi}Q')+
    \Im\kappa\xi  Q' + 
      \Im\frac{\kappa}{\sigma}Q - \omega Q + 
      (1+\Im\delta)|Q|^{2\sigma}Q = 0 \COMMA  \\
&& Q(\xi_1)=\beta \COMMA\;\;\mbox{and} \\
&& Q(\xi)\sim \ASYM{\xi}\;\;\mbox{as}\;\; \xi\to\infty\PERIOD 
\end{eqnarray}
considered in $(\xi_1,\infty)$ has a solution
$$
F(\xi)=F(\beta,\kappa,\omega,\epsilon,\delta,\xi_1;\xi)\PERIOD
$$
Moreover, $F$ can be constructed in such a way that the following conditions
are satisfied
\begin{description}
\item[{\rm (i)}] The mapping from $\LAM_{\rho_1}$ to $\BSPT$ defined by
$$
(\beta,\kappa,\omega,\epsilon,\delta) \mapsto
F(\beta,\kappa,\omega,\epsilon,\delta,\xi_1;\cdot)\;\;\;
$$
is $C^1$ up to the boundary  for each $\vartheta>0$.

\item[{\rm (ii)}] The complex-valued function defined by
$$
(\beta,\kappa,\omega,\epsilon,\delta) \mapsto
\frac{\partial F}{\partial\xi}(\beta,\kappa,\omega,\epsilon,\delta,\xi_1,\xi_1)\;\;\;
$$ 
is $C^1$ (up to the boundary) in $\LAM_{\rho_1}$.
\item[{\rm(iii)}] $F$ and its derivatives $F^{(k)}$ have the following asymptotic expansions
 \begin{eqnarray}
&&F(\xi) =  \ASYM{\xi} 
              \left(\sum_{l=0}^n a_l \xi^{-2l} \, +
                \BIGO(\xi^{-2(n+1)}) \right) \COMMA \label{asymp1}\\
&&F^{(k)}(\xi)  = \frac{\partial^k}{\partial\xi^k}\left(
                 \xi^{-1/\sigma-\Im\omega/\kappa } 
              \sum_{l=0}^n a_l \xi^{-2l} \right)\, +
                \BIGO(\xi^{-1/\sigma-2(n+1)-k})
   \label{asymp2}\PERIOD
\end{eqnarray}
\item[{\rm (iv)}] We have $F(0,\kappa,\omega,\epsilon,\delta,\xi_1,\xi)=0$ and
$$
\frac{\partial F}{\partial\beta}(0,\kappa,\omega,\epsilon,\delta,\xi_1,\xi)=
\frac{P(\kappa,\omega,\epsilon,\xi)}{P(\kappa,\omega,\epsilon,\xi_1)}\PERIOD
$$
\end{description}
\end{theorem}

\noindent
\PROOF See the Appendix \\
\\
\begin{remark}
For $\epsilon>0$ the function $u(x,t)$ given by \VIZ{selfsim} has to be regular at all points $(x,t)$ with $x\neq 0$ by the 
partial regularity theorem proved in \cite{XYan}. Therefore any
solution of \BVPQ\ must admit an asymptotic expansion of the form
stated in (iii) of Theorem~\ref{thm1}. We note that the convergence of
the series $\sum a_l\xi^{-2l}$ in the asymptotic expansion of $Q$ is
equivalent to the analycity of $u$ (in $t$) at the points $(x,T)$,
$x\neq 0$. The asymptotic expansion of the profiles $Q$ does not converge
and therefore $u$ is not analytic in $t$ at any point $(x,T)$.
\end{remark}

Theorem~\ref{thm1} together with elementary perturbation arguments
can be used to show that if  some non-degeneracy conditions are
satisfied, then every solution of \BVPQ\
for $\epsilon=\delta=0$ will persist (with a slight deformation)
for small $\epsilon > 0$ and $\delta>0$.
We will briefly describe this standard procedure for the convenience
of the reader.

First we consider solutions on the finite interval $(0,\xi_1]$ to the initial value problem
\begin{eqnarray}
&& (1-\Im\epsilon)(Q'' + \frac{d-1}{\xi}Q') + \Im\kappa\xi Q'
 + \Im\frac{\kappa}{\sigma} Q - \omega Q + (1+\Im\delta)|Q|^{2\sigma}Q
 = 0  \COMMA \label{IVP}\\
&& Q(0) = \mu \COMMA\;\;\; Q'(0)=0 \PERIOD \label{IVPa}
\end{eqnarray}

We denote the solution (if it exists) by $G(\xi) =
G(\mu,\NU,\ALPHA;\xi)$. Clearly the set of parameters $(\mu,\NU,\ALPHA)$ for
which $G$ is well-defined is open. We define
\begin{equation}\label{clg15}
\beta(\mu,\NU,\ALPHA) = G(\mu,\NU,\ALPHA;\xi_1)\PERIOD
\end{equation}

Assume that the boundary-value problem \BVPQ\  has a solution 
which satisfies $Q(\xi_1)$$=\mu$. 
With a slight abuse of notation, let us denote such a solution
by  $Q(\mu,\NU,\ALPHA;\xi)$. 
Clearly $Q(\mu,\NU,\ALPHA;\xi)$ is defined only on a submanifold of the
parameter space, but this will not
be important in what follows.
As we have seen in the introduction, we have
$$
\lambda^{1/\sigma+\Im\theta} Q(\mu,\NU,\ALPHA,\lambda\xi) = 
Q(\lambda^{1/\sigma+\Im\theta}\mu,\lambda^2\kappa,\lambda^2\omega,\ALPHA;\xi)\COMMA
$$
for all $\lambda>0$ and $\theta\in\PER$. Therefore we can work
with normalized solutions, i.\ e.\ we assume
that $\omega=1$  and that  $Q(0)$ is real and non-negative. 
Assume $\mu_0>0, \kappa_0>0$ and
suppose $Q(\mu_0,\kappa_0,1,0,0;\xi)$ exists.

We set
\begin{equation}\label{clg16}
g(\mu,\kappa,\ALPHA) = \PD{G}{\xi}(\mu,\kappa,1,\ALPHA;\xi_1) - 
                     \PD{F}{\xi}(\beta(\mu,\kappa,1,\ALPHA),\kappa,1,\ALPHA,\xi_1,\xi_1) 
\PERIOD
\end{equation}
where $\beta(\mu,\NU,\ALPHA)$ is defined by \VIZ{clg15}.
By our assumptions and by Theorem~\ref{thm1} the mapping $g$ is
well defined and continuously differentiable in
a set of the form $(\mu_1,\mu_2)\times(\omega_1,\omega_2)\times
[0,\epsilon_1)\times(\delta_1,\delta_2)$ containing 
the point $(\mu_0,\kappa_0,0,0)$.
Since
$g(\mu_0,\kappa_0,0,0)=0$, we see that the equation 
$$ 
g(\mu,\kappa,\epsilon,\delta)=0
$$
has solutions for small $\epsilon>0,\,\delta>0$ if the following
non-degeneracy condition is satisfied:

\begin{equation}\label{clgC}
\DET\left(\begin{array}{cc}
                    g^1_\mu  & g^1_\kappa \\
                    g^2_\mu  & g^2_\kappa
                 \end{array} \right)
    \neq 0 \SPACE\mbox{at}\SPACE (\mu_0,\kappa_0,0,0)\COMMA
\end{equation} 
where $g^1 = \REAL{g}$, $g^2=\IMAG{g}$ and subscripts denote partial
derivatives with respect to the corresponding variables.

Thus a non-trivial solution of \BVPQ\  for $\epsilon=\delta=0$ also
gives a solution of \BVPQ\  for $\epsilon,\delta>0$ if \VIZ{clgC}
is satisfied. Based on our numerical calculations described in the
next section, we conjecture that \VIZ{clgC} is satisfied
for $\sigma>2/d$.

\section{Numerical results}\label{sec3}
Theorem~\ref{thm1} allows us to rewrite the boundary value problem
\BVPQ\ as a boundary value problem on a finite interval $(0,\xi_1)$
in the following way.
\begin{eqnarray}
&& (1-\Im\epsilon)(Q'' + \frac{d-1}{\xi} Q') +
      \Im\kappa\xi Q' + 
      \Im\frac{\kappa}{\sigma} Q -  Q + 
      (1+\Im\delta)|Q|^{2\sigma}Q = 0 \COMMA \label{clg17}\\
&& Q'(0)=0 \COMMA \\
&& Q(\xi_1) = \beta \COMMA\\
&& Q'(\xi_1) = \PD{F}{\xi}(\beta,\kappa,1,\ALPHA,\xi_1;\xi_1) \COMMA \label{clg17a}
\end{eqnarray}
where the unknown quantities are $Q$, $\beta$ and $\kappa$. 
Of course, the problem \VIZ{clg17}-\VIZ{clg17a} is equivalent 
to the equation $g=0$ in the  previous section.
If we approximate $F$ by the first few terms of its asymptotic expansion, 
the problem \VIZ{clg17}-\VIZ{clg17a} can be solved numerically. 
In our numerical computations we
investigated the dependence on $\xi_1$ and on the number of terms of the
asymptotic expansion of $F$. It turned out that $\xi_1\sim 30$ and
the first term of the asymptotic expansion already worked very well.
However, many of our computations were done with the first two terms of the 
asymptotic expansion \VIZ{asymp1}. 
In the case $\epsilon=\delta=0$, i.e. NLS, the 
value of $(\mu_1,\kappa_1)$ was computed in 
\cite{LeMesurierPapanicolau} using
completely different approach. The values presented in that paper are
in an excellent agreement with our computations, see below.

The first term of the asymptotic expansion of $F$ is
$F\sim \beta(\xi/\xi_1)^{-1/\sigma-\Im\omega/\kappa}$. Using this
approximation, we obtain from \VIZ{clg17}-\VIZ{clg17a}
\begin{eqnarray}
&& (1-\Im\epsilon)(Q'' + \frac{d-1}{\xi} Q') +
      \Im\kappa\xi Q' + 
      \Im\frac{\kappa}{\sigma} Q -  Q + 
      (1+\Im\delta)|Q|^{2\sigma}Q = 0 \COMMA \label{clg18}\\
&& Q'(0)=0 \COMMA \label{clg18a} \\
&& \xi_1 Q'(\xi_1) + \left(\frac{1}{\sigma} +
  \Im\frac{1}{\kappa}\right) Q(\xi_1) = 0 \PERIOD \label{clg18b}
\end{eqnarray}
Higher order approximations can be derived in a similar way.
 Note that in the formulation of \VIZ{clg18} we already
fixed the normalization $\omega=1$, so that the unknowns are $Q$
and $\kappa$. The boundary condition \VIZ{clg18b} is also closely related to the
boundary condition used in \cite{Zacharov2}, \cite{LandmanPapanicolaou} for
simulations based on solving time dependent problem in the PDE
\VIZ{clg1}-\VIZ{clg1b}. 

There are essentially two approaches to the numerical solution of the
boundary-value problem \BVPS. One can use collocation methods to approximate 
the boundary-value problem and then to apply Newton's method to the
discretization of the non-linear operator which defines the equation
\VIZ{clg18}. Implementation of this strategy requires further changes
in the formulation since the linearization of the non-linear operator
always has zero in its spectrum due to the $S^1$-equivariance of the
equation. Therefore a shooting method was easier to
implement and it  also proved to be sufficiently accurate. Because of 
well-known
sensitivity of shooting methods to problem parameters we performed
computations in different normalizations: with $Q(0)=1$ fixed
and parameters $(\kappa,\omega)$ as the unknowns as well as 
with  $\omega=1$ fixed and parameters
$(\kappa,\mu)$ as the unknowns. Moreover we compared both backward 
and forward shooting methods on the interval $(0,\xi_1)$. 
All computed solutions turned out to be in a very good agreement.

Here we describe only the shooting method using the normalization
$\omega=1$, in which we calculate $\xi_1 Q'(\xi_1)+ (1/\sigma +
\Im/\kappa) Q(\xi_1)$ as a function of the parameters $\mu=Q(0),
\kappa,\epsilon$, and $\delta$,
i.e., we solve the equation
\begin{equation}
f(\mu,\kappa,\epsilon,\delta) \equiv 
\xi_1 Q'(\xi_1)+ \left(1/\sigma +\Im/\kappa\right) Q(\xi_1) = 0 \PERIOD
\end{equation}
The integration of the underlying ODE must be done with
sufficient accuracy. We compared various ODE solvers. 
A variable-order, variable-step Adams method 
as implemented, for example, in the NAG
library proved to be sufficiently accurate in most of the
calculations. 
To locate the initial values for each branch we inspected the
two dimensional subspace of the parameter space given by $(\mu,\kappa,0,0)$ 
for $d=1$ and $\sigma=2.3$
and computed the degree of the
function $f$ restricted to that subspace along various curves.
The other solutions were then calculated by continuation.
This was done with the help of bifurcation analysis 
package developed as a part of \cite{Ple}. 
The implementation of the path-following
procedure with a Newton corrector step
can be done efficiently as the linearization along a solution
is evaluated at the same integration step as the solution.

As we described in the introduction, it appears that the equation
$f(\mu,\kappa,0,0)$ has countably many solutions $(\kappa_j,\mu_j)$ in the
region $\kappa>0$, $\mu>0$. The corresponding profiles $|Q_j|$, when extended
to $(-\infty, \infty)$ as even functions, have exactly $j$ local maxima.
\begin{remark}
As mentioned in the beginning of this section the solutions
$(\kappa_j,\mu_j)$ may slightly depend on the value $\xi_1$ 
and the approximation of $F$.  For  
solutions  described here we tested the dependence of the results
on $\xi_1$ for  $\xi_1\in[20,100]$, and also on the approximation of $F$
by taking either one or two terms in the asymptotics expansion of $F$.
We also directly compared our numerical solutions in intervals of the form
$(\eta_1,\xi_1)$
with the explicit formulae given by one or two terms of the asymptotic
expansion of the solution in $(\eta_1,\xi_1)$ for various values of $\eta_1$.
All these tests indicated a good convergence of our approximations.
Based on these tests, we estimate that the error in the values of 
the ``roots'' $(\kappa_j,\mu_j)$ is of the order $10^{-3}$ or better.
\end{remark}

\noindent
{\bf Case I} ($d=1$, $\sigma=2.3$, $\delta=0$): Results for the 
one-dimensional case are tabulated in Table~\ref{tab1}. 

\begin{table}[ht]
\begin{center}
\begin{tabular}{|c|c|ll|ll|}\hline
branch       &  turning point  & 
\multicolumn{2}{c|}{$\omega=1$}& \multicolumn{2}{c|}{$Q(0)=1$} \\ \cline{3-6}
\multicolumn{1}{|c|}{$j$}      & \multicolumn{1}{c|}{$\epsilon^*$}             &
\multicolumn{1}{c}{$\kappa $}  & \multicolumn{1}{c|}{$\mu$} &    
\multicolumn{1}{c}{$\kappa $}  & \multicolumn{1}{c|}{$\omega$}\\ 
\hline
\num{1}& \num{0.06064}& \num{0.85311}& \num{1.23204}& \num{0.32669}& \num{0.38294}\\
\num{2}& \num{0.05182}& \num{0.49323}& \num{0.78308}& \num{1.51894}& \num{3.07959}\\
\num{3}& \num{0.04466}& \num{0.34673}& \num{1.12388}& \num{0.20263}& \num{0.58438}\\
\num{4}& \num{0.03900}& \num{0.26678}& \num{0.78308}& \num{0.47127}& \num{1.76651}\\
\num{5}& \num{0.03455}& \num{0.21643}& \num{1.07947}& \num{0.15225}& \num{0.70345} \\
\num{6}& \num{0.03099}& \num{0.18185}& \num{0.92714}& \num{0.25750}& \num{1.41624} \\
\num{7}& \num{0.02803}& \num{0.15667}& \num{1.05430}& \num{0.12284}& \num{0.78409}\\
\num{8}& \num{0.02559}& \num{0.13756}& \num{0.95061}& \num{0.17365}& \num{1.26236}\\
\hline
\end{tabular}
\end{center}
\caption{{\sc Case I} ($d=1$, $\sigma=2.3$, $\delta=0$):
 $\epsilon^*$ defines the position of the turning point; coordinates
 $(\kappa,\mu)$, where $\mu=Q(0)$ define the initial point of each branch, when
$\epsilon=0$, in the normalization $\omega=1$.
 The coordinates $(\kappa,\omega)$ refer to the 
same solutions in the  normalization $Q(0)=1$.}\label{tab1}
\end{table}

\begin{figure}[ht]
\begin{center}
\begin{minipage}{0.48\linewidth}
\centering\epsfig{file=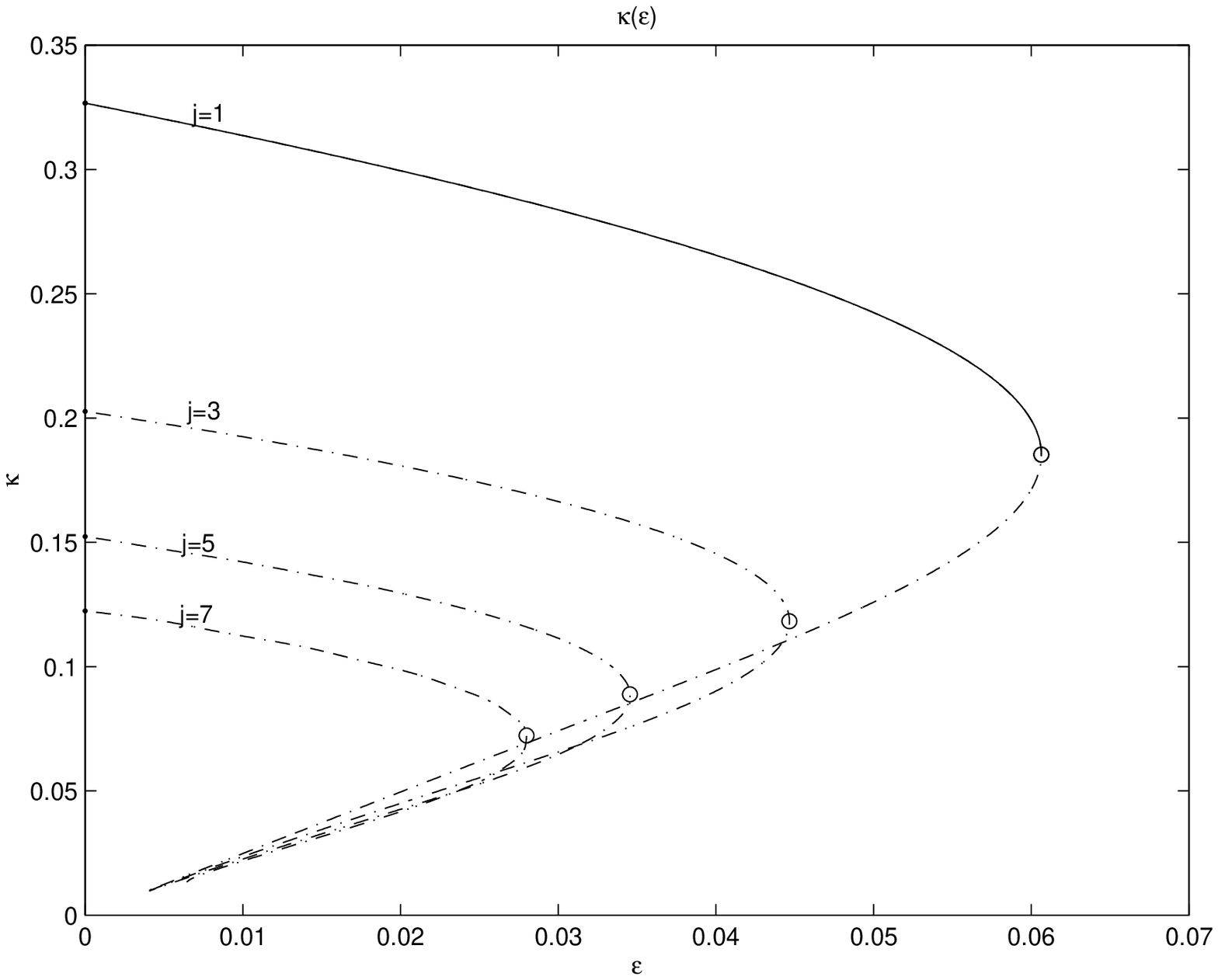, width=.98\linewidth}
\end{minipage}
\begin{minipage}{0.48\linewidth}
\centering\epsfig{file=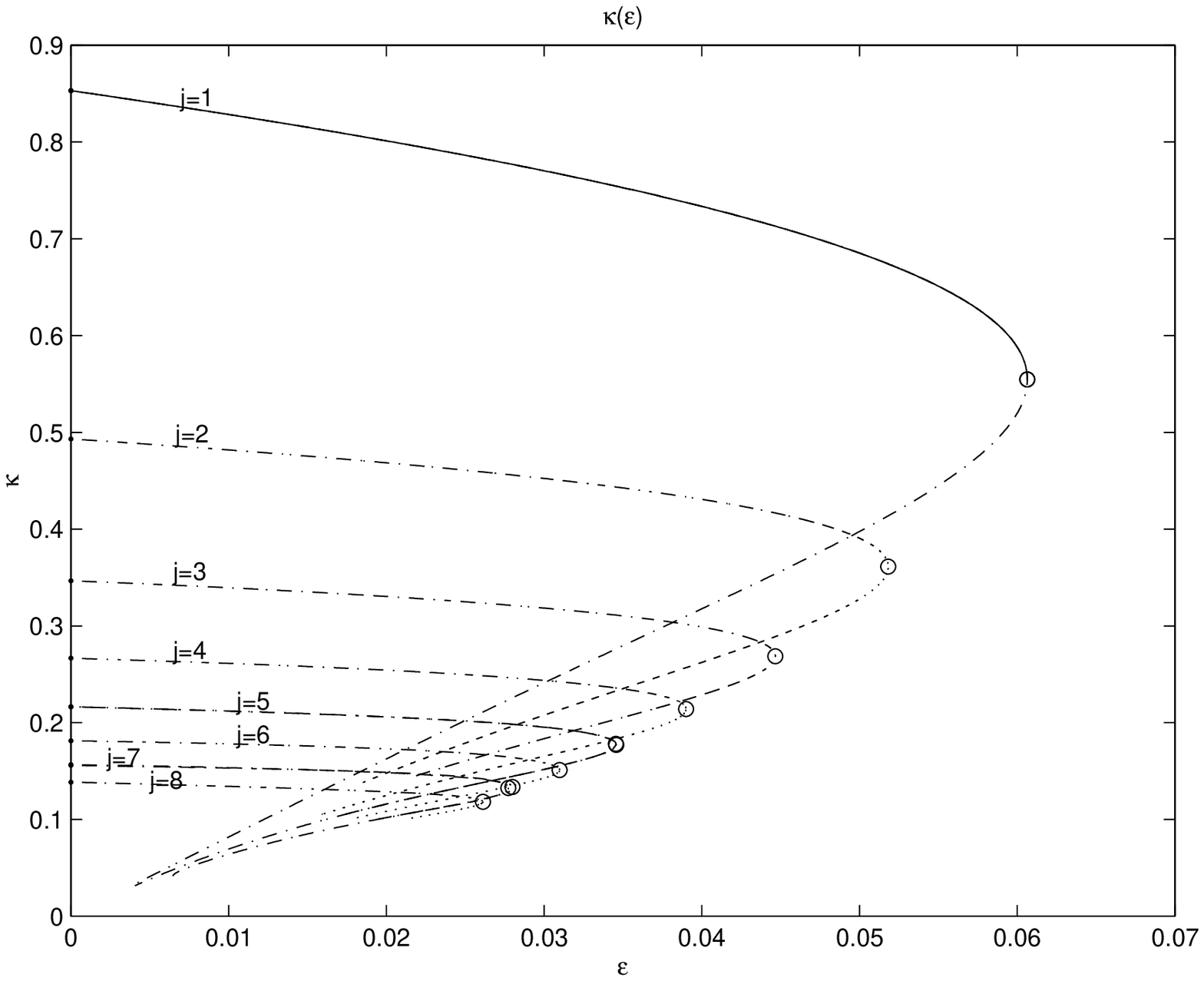, width=.98\linewidth}
\end{minipage}
\end{center}
\caption{{\sc Case I} ($d=1$, $\sigma=2.3$, $\delta=0$): solution branches
  $(\epsilon,\kappa_j(\epsilon))$; normalization $Q(0)=1$ (left) and
  $\omega=1$ (right). Turning points are denoted
  by $o$. The solid line indicates stable solutions, the dashed line
  indicates unstable solutions. See Section~\ref{sec4} for details
concerning stability issues.}\label{fig1}
\end{figure}

\begin{figure}[ht]
\begin{center}
\begin{minipage}{0.48\linewidth}
\centering\epsfig{file=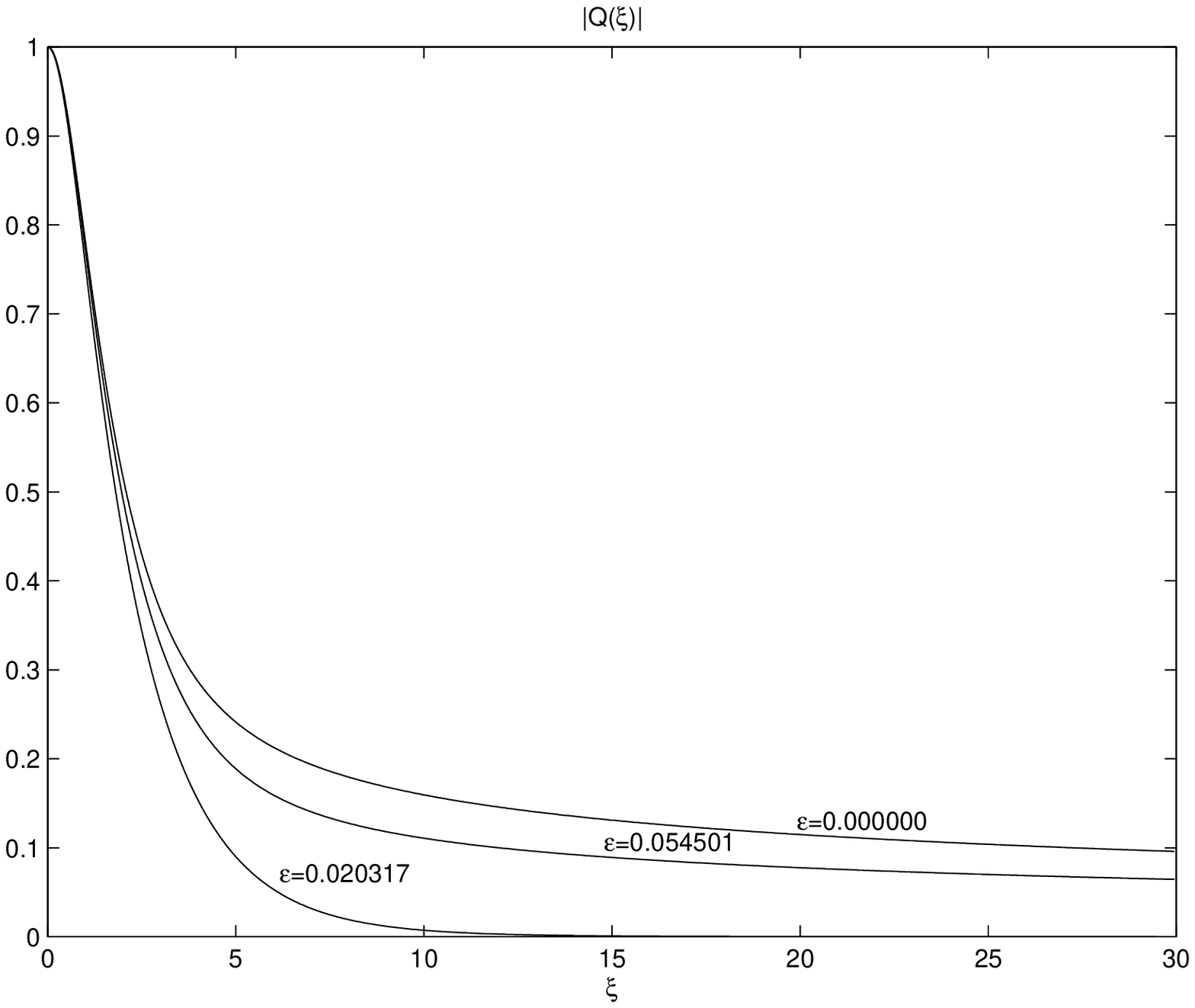, width=.98\linewidth}
\end{minipage}
\begin{minipage}{0.48\linewidth}
\centering\epsfig{file=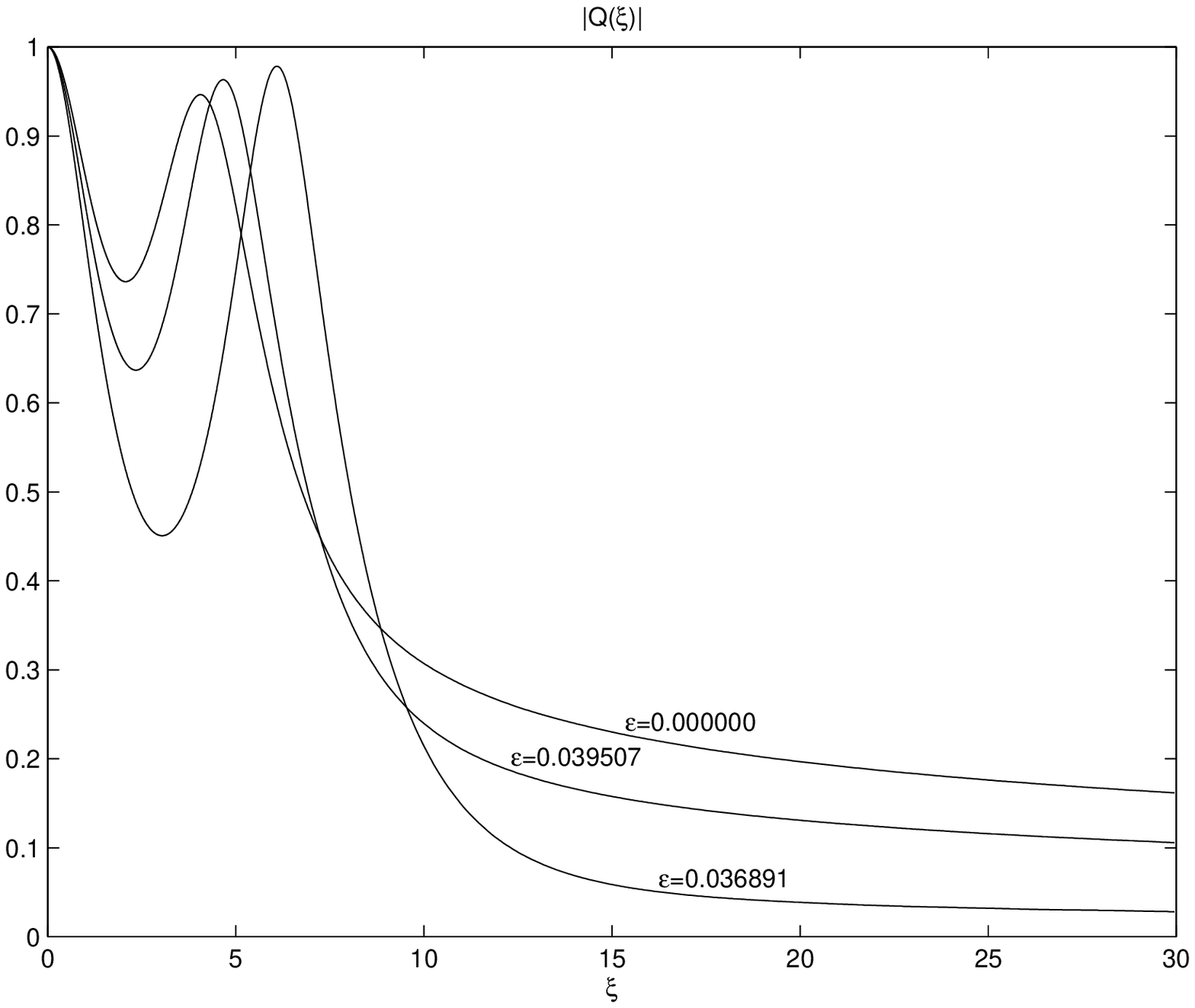, width=.98\linewidth}
\end{minipage}
\end{center}
\caption{{\sc Case I} ($d=1$, $\sigma=2.3$, $\delta=0$):  profiles $|Q(\xi)|$ at points along
         the $j$-th branch; $j=1$ (left), $j=3$ (right); normalization $Q(0)=1$.}\label{figd1-b13}
\end{figure}

\begin{figure}[ht]
\begin{center}
\begin{minipage}{0.48\linewidth}
\centering\epsfig{file=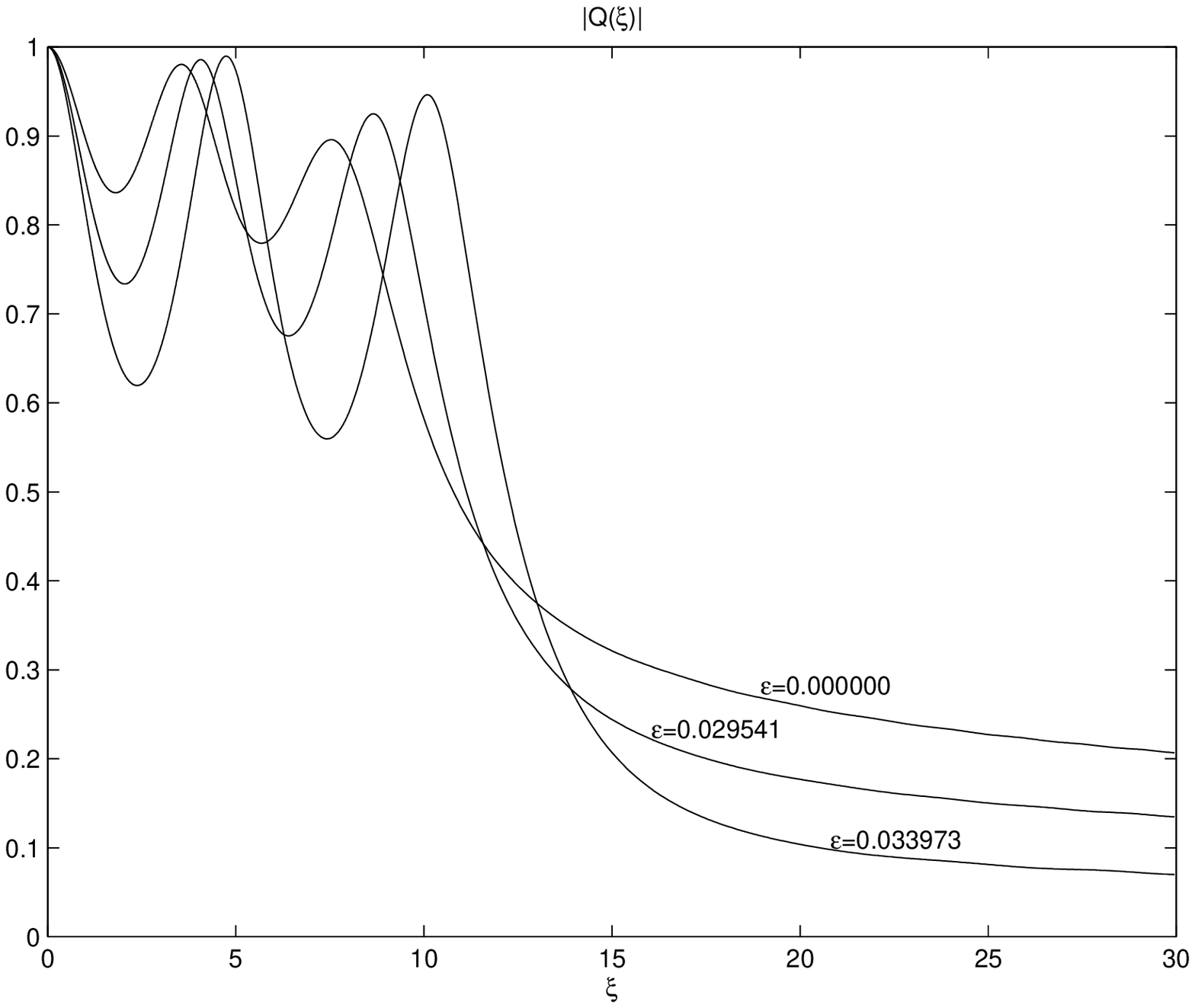, width=.98\linewidth}
\end{minipage}
\begin{minipage}{0.48\linewidth}
\centering\epsfig{file=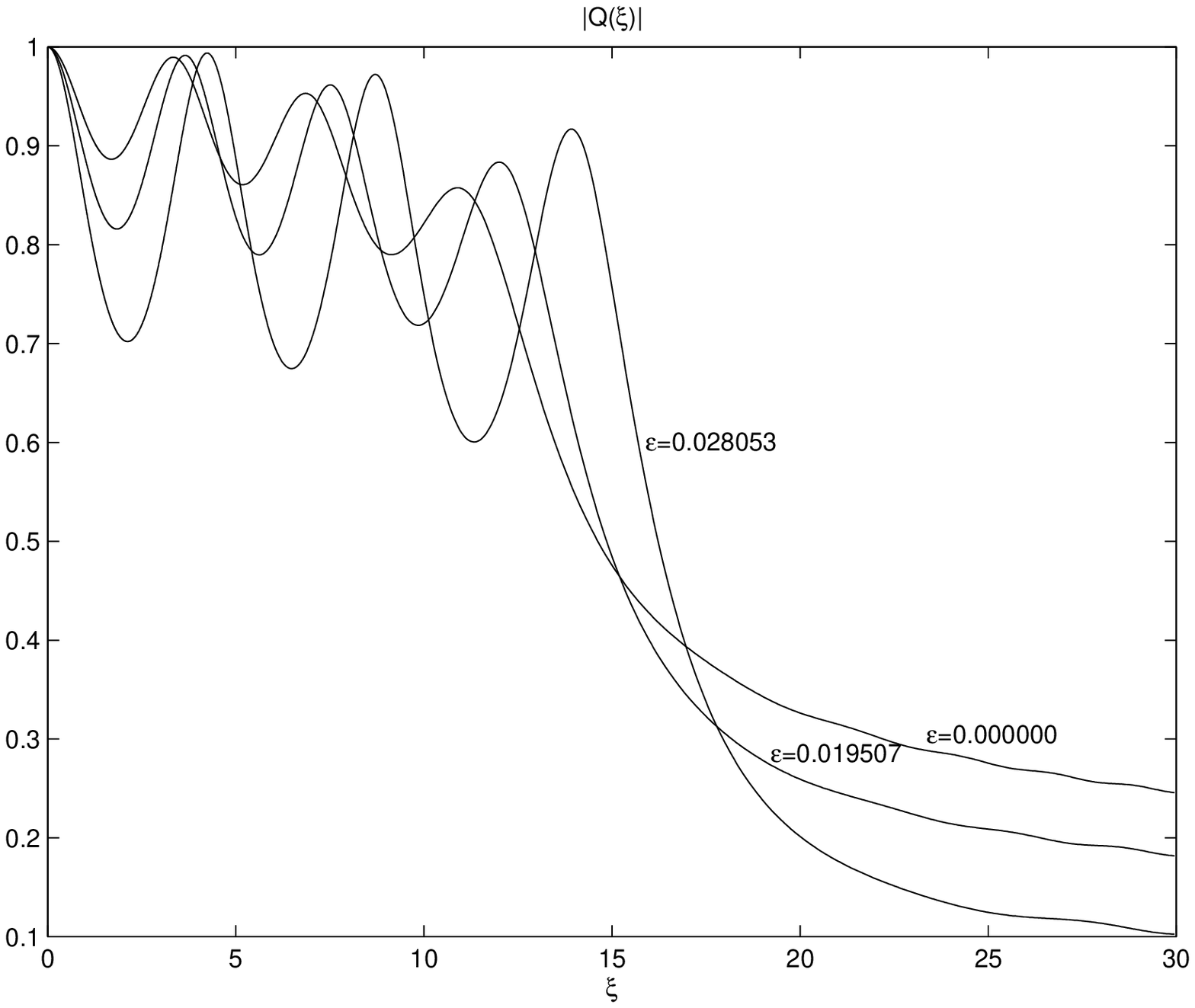, width=.98\linewidth}
\end{minipage}
\end{center}
\caption{{\sc Case I} ($d=1$, $\sigma=2.3$, $\delta=0$): profiles $|Q(\xi)|$ at points along
         the $j$-th branch; $j=5$ (left), $j=7$ (right); normalization $Q(0)=1$.}\label{figd1-b57}
\end{figure}

\begin{figure}[ht]
\begin{center}
\begin{minipage}{0.48\linewidth}
\centering\epsfig{file=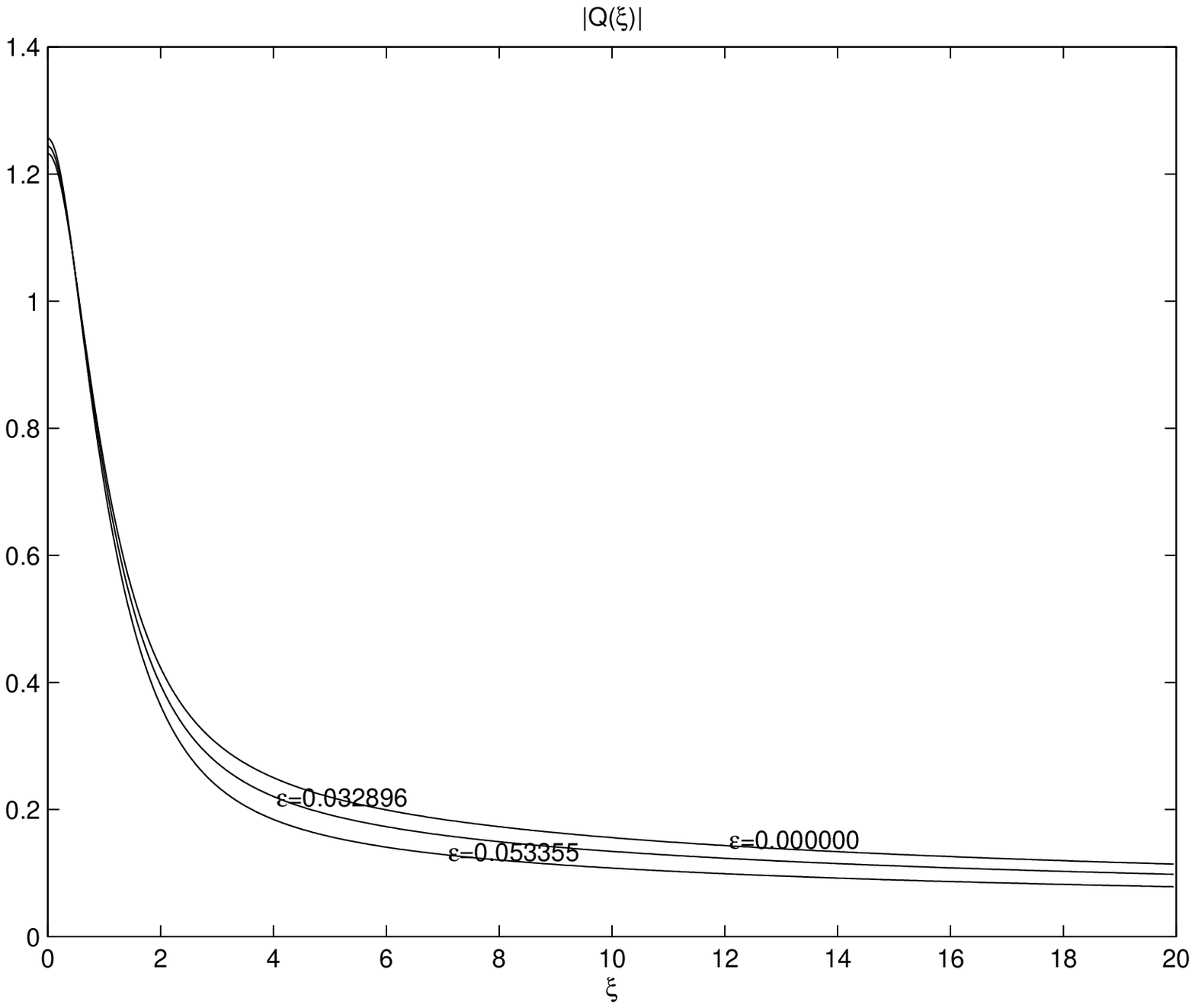, width=.98\linewidth}
\end{minipage}
\begin{minipage}{0.48\linewidth}
\centering\epsfig{file=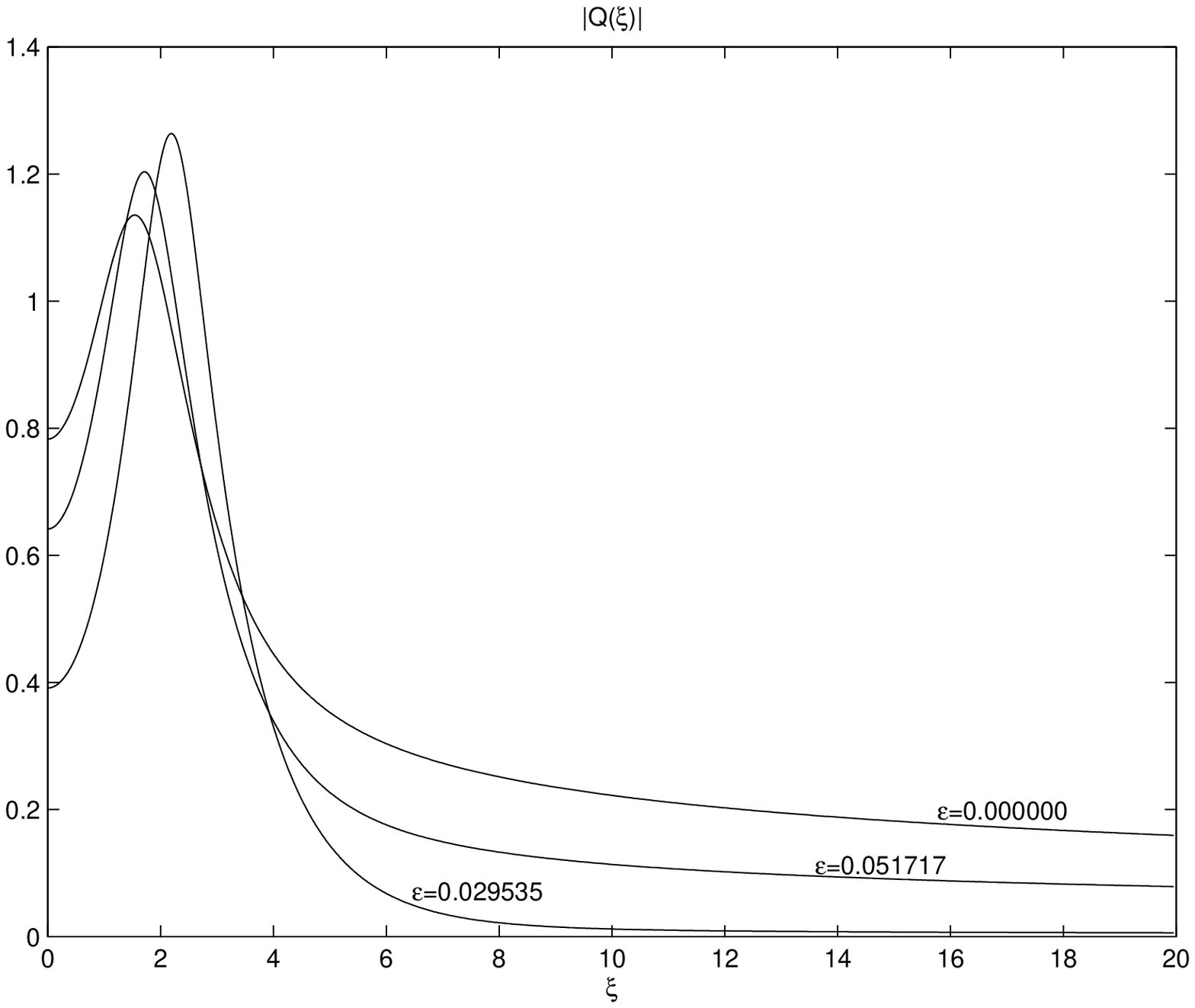, width=.98\linewidth}
\end{minipage}
\end{center}
\caption{{\sc Case I} ($d=1$, $\sigma=2.3$, $\delta=0$): profiles $|Q(\xi)|$ at points along
         the $j$-th branch; $j=1$ (left), $j=2$ (right); normalization $\omega=1$.}
        \label{figd1-bmu13}
\end{figure}

\begin{figure}[ht]
\begin{center}
\begin{minipage}{0.48\linewidth}
\centering\epsfig{file=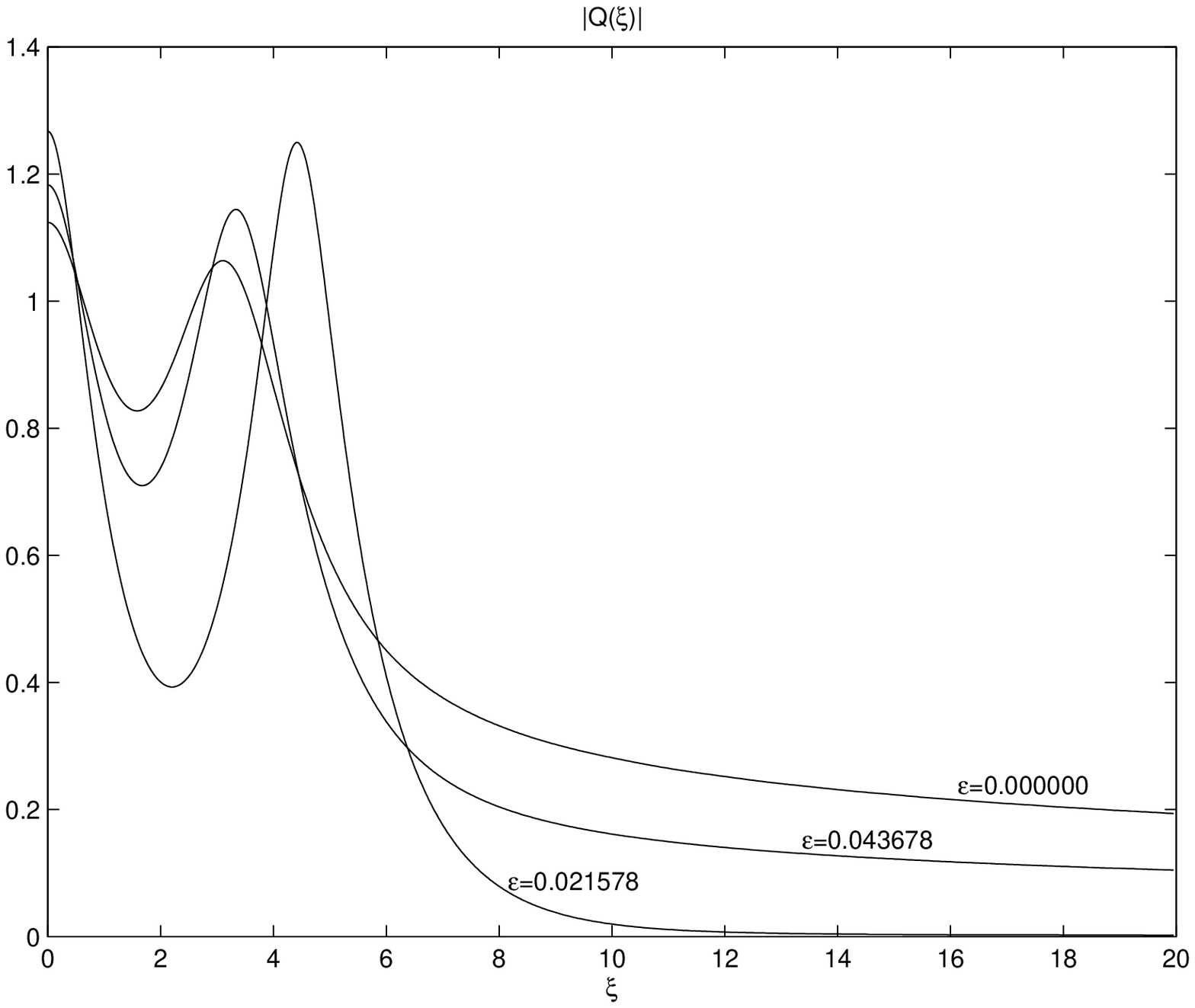, width=.98\linewidth}
\end{minipage}
\begin{minipage}{0.48\linewidth}
\centering\epsfig{file=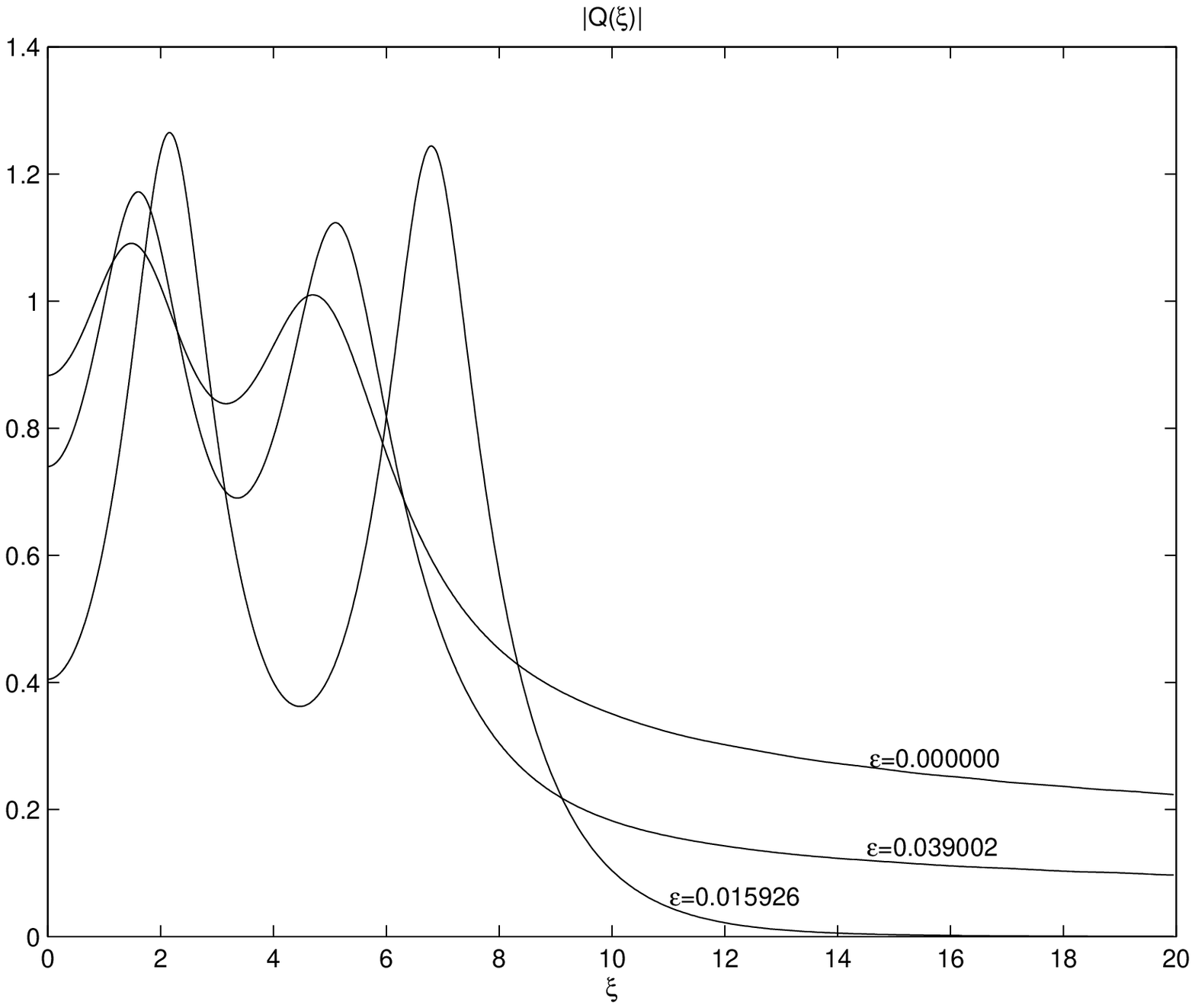, width=.98\linewidth}
\end{minipage}
\end{center}
\caption{{\sc Case I} ($d=1$, $\sigma=2.3$, $\delta=0$): profiles $|Q(\xi)|$ at points along
         the $j$-th branch; $j=3$ (left), $j=4$ (right); normalization $\omega=1$.}\label{figd1-bmu57}
\end{figure}

The continuation
of solutions parameterized by $\epsilon$ is depicted in
Figure~\ref{fig1}. To give the reader a good idea about the form of
the solutions we plot profiles $|Q(\xi)|$ at a few points on each
branch
(at the point $\epsilon=0$ (solution to NLS), at a point on the upper
part of the branch and at another point on the lower part of the
branch). 
We used both normalizations, $Q(0)=1$ in
Figure~\ref{figd1-b13}--\ref{figd1-b57}, and $\omega=1$ in 
Figure~\ref{figd1-bmu13}--\ref{figd1-bmu57}. Normalization $Q(0)=1$ is
convenient for computing solutions along the odd branches as these
solutions
exhibit a maximum at the origin which is also present in the solution 
for $\epsilon\to 0$ on the lower part of the branch.

We now describe the behavior of solutions along the branches.
 We observe that as $\epsilon$
returns to zero after passing through the turning point, 
we have  $\kappa_j(\epsilon)\to 0$. 

For $j=1$ 
the solution $Q^\epsilon_1(\xi)$ approaches a specific solution, 
usually called {\it the ground
  state}, of the equation
\begin{equation}\label{clg22}
- u'' + u = |u|^{2\sigma}u \PERIOD
\end{equation}
For $j>1$ the profile $|Q^\epsilon_j(\xi)|$ seems to separate into $j$ copies of the
ground-state solution which are moving to infinity. For $j$ even
all of them ``escape'' to infinity, while for $j$ odd one will stay at 
the origin and the rest will move to the infinity as $\epsilon\to 0$ and
$\kappa_j(\epsilon)\to 0$.

\noindent
{\bf Case II} ($d=3$, $\sigma=1$, $\delta=0$): Results for the
three-dimensional  case are tabulated in Table~\ref{tab2}.

\begin{table}[ht]
\begin{center}
\begin{tabular}{|c|c|ll|ll|}\hline
branch       & turning point   & 
\multicolumn{2}{c|}{$\omega=1$}& \multicolumn{2}{c|}{$Q(0)=1$} \\ \cline{3-6}
\multicolumn{1}{|c|}{$j$}      & \multicolumn{1}{c|}{$\epsilon^*$}             &
\multicolumn{1}{c}{$\kappa $}  & \multicolumn{1}{c|}{$\mu$} &    
\multicolumn{1}{c}{$\kappa $}  & \multicolumn{1}{c|}{$\omega$}\\ 
\hline
\num{1}& \num{0.19813}& \num{0.91737}& \num{1.88529}& \num{0.25810}& \num{0.28135}\\
\num{2}& \num{0.24402}& \num{0.32091}& \num{0.83559}& \num{0.45535}& \num{1.41727}\\
\num{3}& \num{0.22762}& \num{0.22704}& \num{1.10834}& \num{0.18242}& \num{0.80684}\\
\num{4}& \num{0.19520}& \num{0.16543}& \num{1.03257}& \num{0.15516}& \num{0.93792}\\
\num{5}& \num{0.18168}& \num{0.14237}& \num{1.00325}& \num{0.13241}& \num{0.96677}\\
\hline
\end{tabular}
\end{center}
\caption{{\sc Case II}  ($d=3$, $\sigma=1$, $\delta=0$):
$\epsilon^*$ defines the position of the turning point;
coordinates $(\kappa,\mu)$, where $\mu=Q(0)$ define the 
initial point of each branch, when $\epsilon=0$, 
in the normalization $\omega=1$.
The coordinates $(\kappa,\omega)$ refer to the same solutions 
in the  normalization $Q(0)=1$. 
In \cite{LeMesurierPapanicolau} the solution corresponding to $j=1$ 
and $\epsilon=0$ was calculated by the method of dynamical 
rescaling, which gave $\kappa=0.917$ and $\mu=1.885$.} \label{tab2}
\end{table}

Some data for NLS ($\epsilon=\delta=0$) and the ``basic solution''
(corresponding to the beginning of our first branch) are available in the
literature and can be used to estimate accuracy of our calculations.
One can see a very good agreement of our solution for $j=1$ with values  
$\kappa=\num{0.917}$ and
$\mu\equiv Q(0)=\num{1.885}$ obtained in
\cite{LeMesurierPapanicolau} from simulations that used the dynamical rescaling
method applied to the initial value problem \VIZ{clg1a}-\VIZ{clg1b}. 

\begin{figure}[ht]
\begin{center}
\begin{minipage}{0.48\linewidth}
\centering\epsfig{file=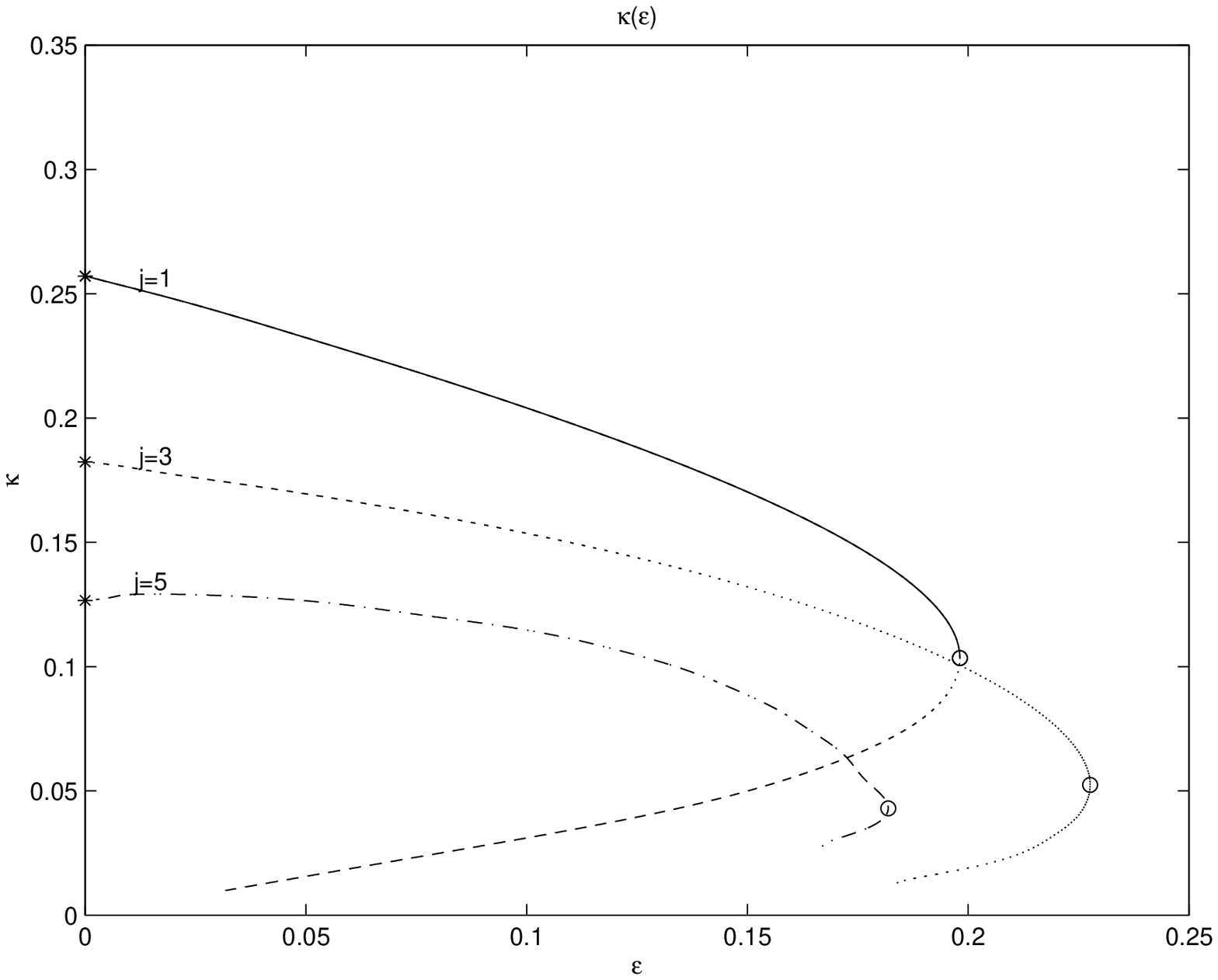, width=.98\linewidth}
\end{minipage}
\begin{minipage}{0.48\linewidth}
\centering\epsfig{file=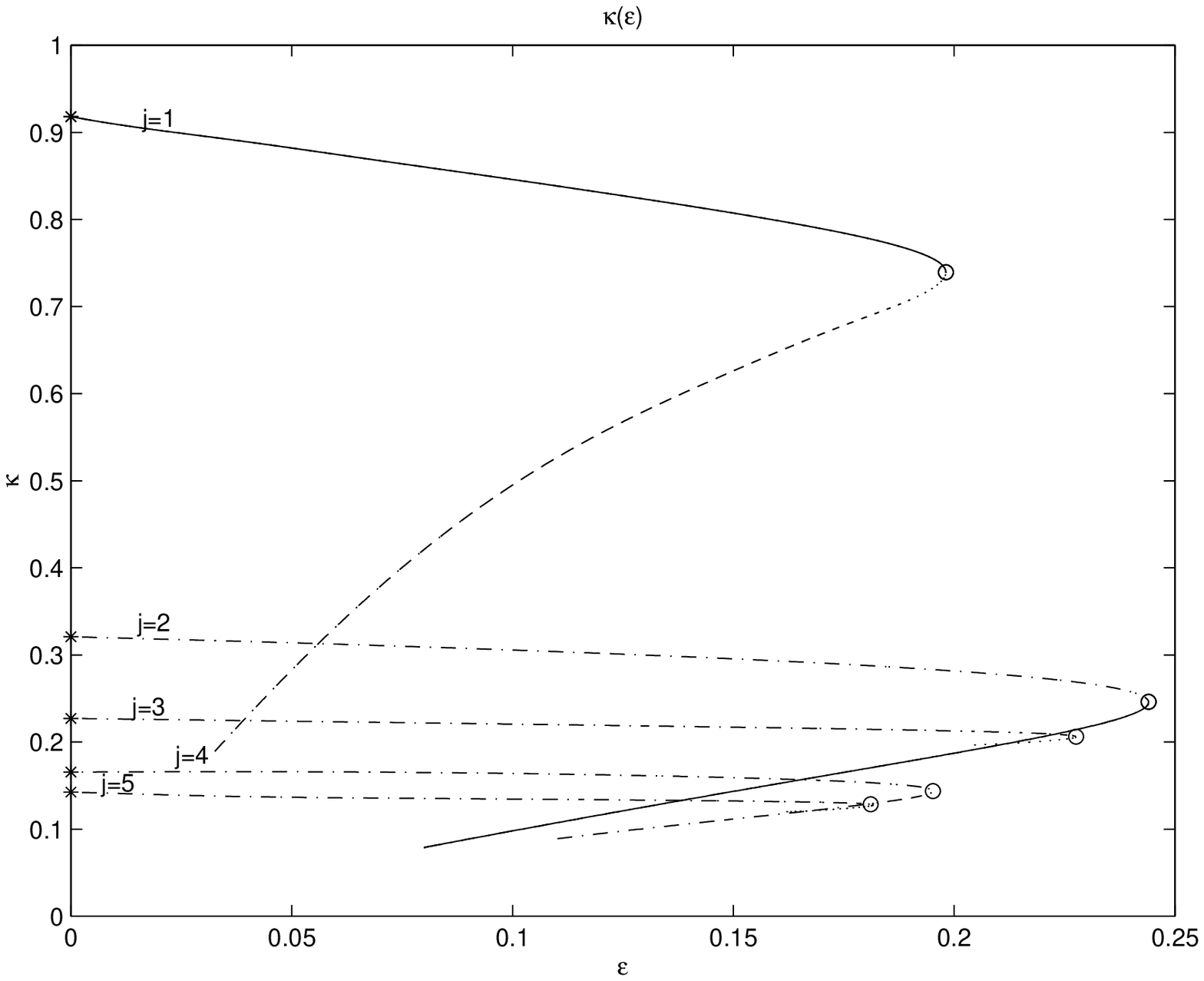, width=.98\linewidth}
\end{minipage}
\end{center}
\caption{{\sc Case II} ($d=3$, $\sigma=1$, $\delta=0$): solution branches
  $(\epsilon,\kappa_j(\epsilon))$; normalization $Q(0)=1$ (left) and
  $\omega=1$ (right). Turning points are denoted
  denoted by $o$. The solid line indicates stable solutions, the dashed line
  indicates unstable solutions. See Section~\ref{sec4} for details regarding
  stability issues.}\label{fig2}
\end{figure}

The behavior of solutions $Q^\epsilon_j$ when $\epsilon$ returns to zero is
qualitatively similar to the one-dimensional case, but there are also
some new interesting features.

First, the maximal possible $\epsilon^*_j$ is not attained on the first branch
($j=1$).
Second, the behavior as $\epsilon\to 0$ and $\kappa_j(\epsilon)\to 0$ 
is different for even and odd values $j$. For $j=2k-1$ (odd) the solutions
$Q^\epsilon_j(\xi)$ converge (along the $j$th-branch) to the $k$-th radial
solution of the problem
\begin{eqnarray}\label{clg23}
-\LAP u + u & = & |u|^{2\sigma}u \;\;\mbox{in}\;\;\RR{3} \COMMA\\
          u & \to & 0 \;\;\mbox{as}\;\;|x|\to \infty \PERIOD
\end{eqnarray}

For $j=2k$ (even), one observes a similar behavior as in the
one-dimensional case.
\begin{figure}[ht]
\begin{center}
\begin{minipage}{0.48\linewidth}
\centering\epsfig{file=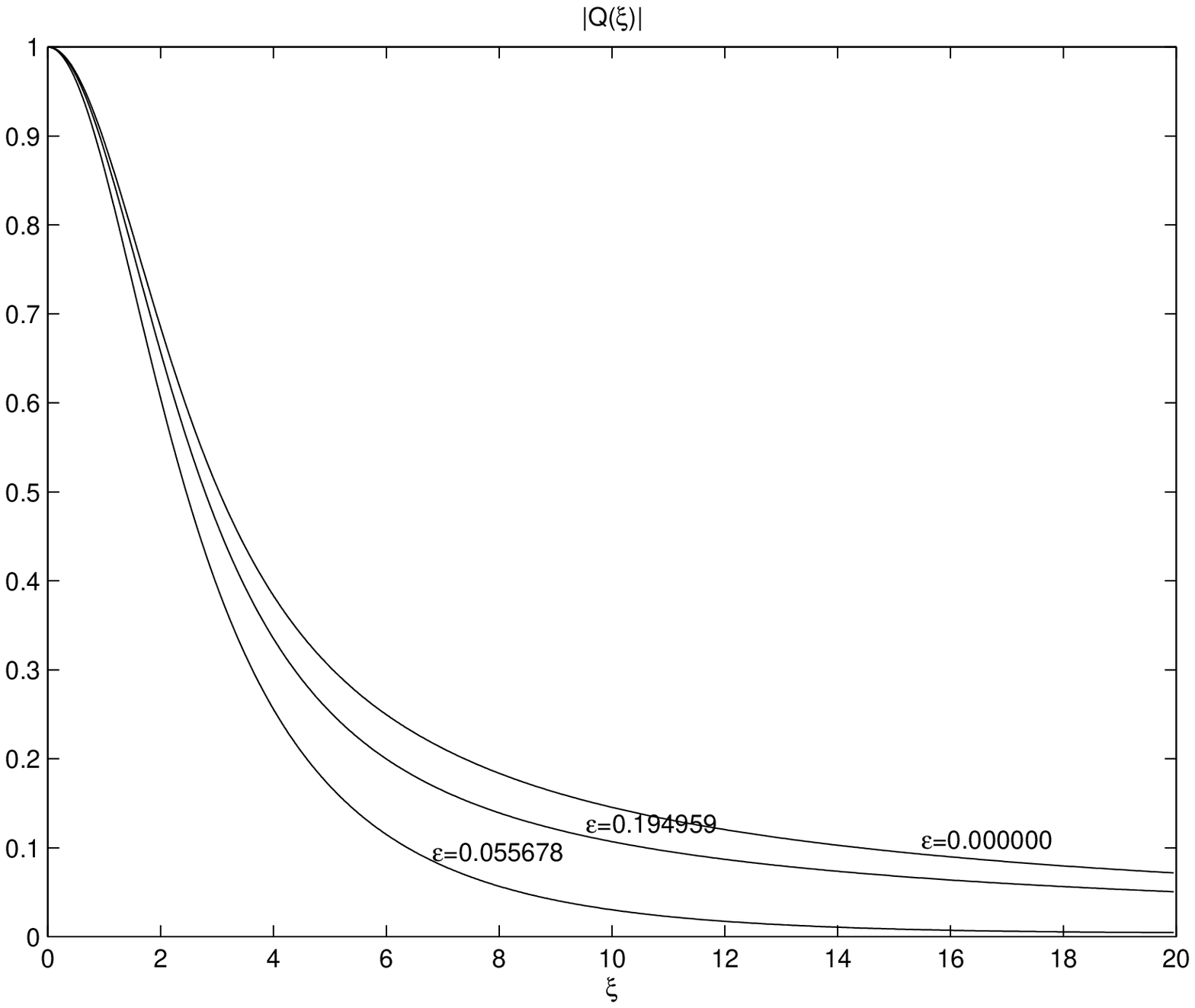, width=.98\linewidth}
\end{minipage}
\begin{minipage}{0.48\linewidth}
\centering\epsfig{file=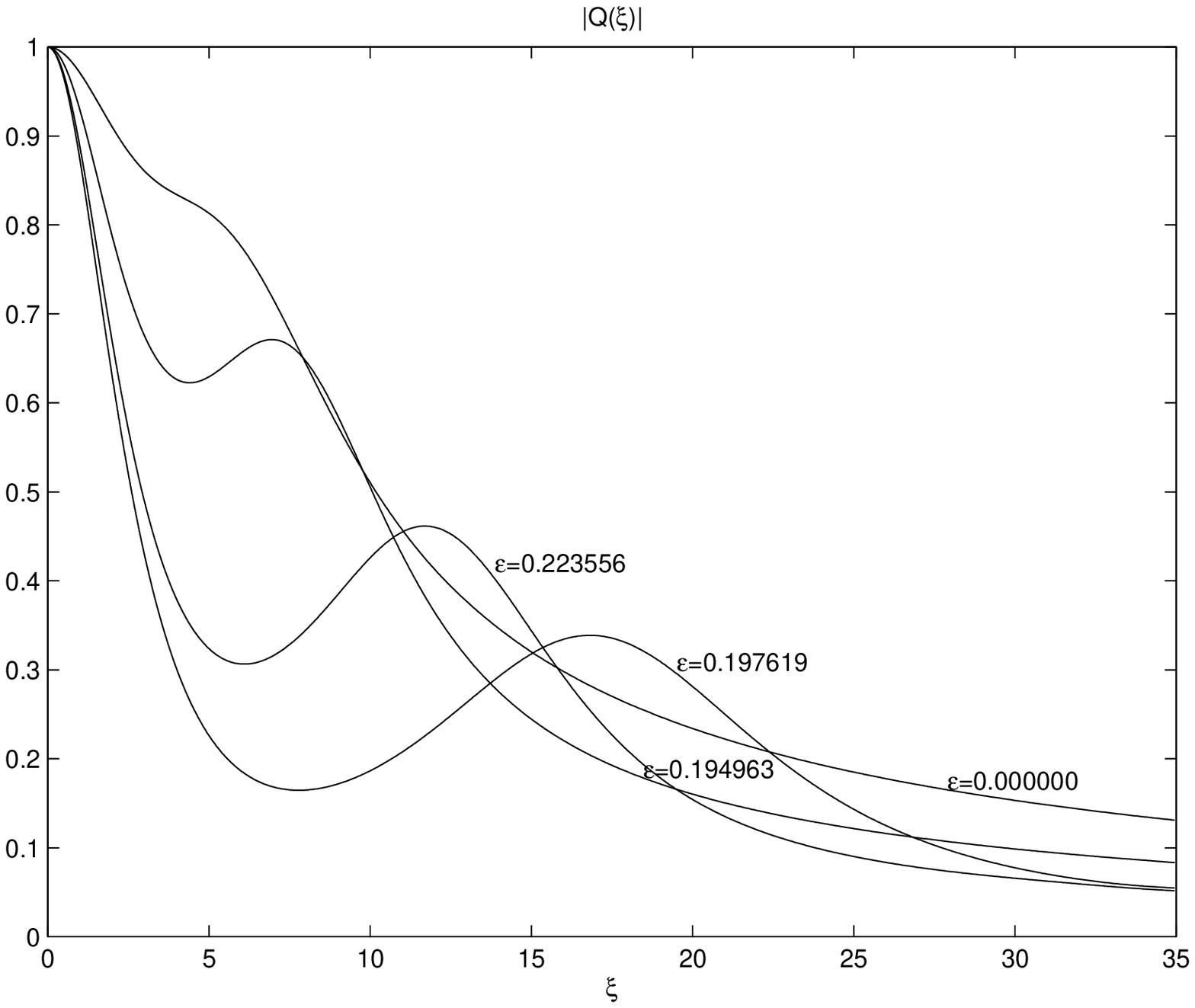, width=.98\linewidth}
\end{minipage}
\end{center}
\caption{{\sc Case II} ($d=3$, $\sigma=1$, $\delta=0$): the profile $|Q(\xi)|$ at points along
         the $j$-th branch; $j=1$ (left), $j=3$ (right); normalization $Q(0)=1$.}\label{figd3-b13}
\end{figure}

\begin{figure}[ht]
\begin{center}
\begin{minipage}{0.48\linewidth}
\centering\epsfig{file=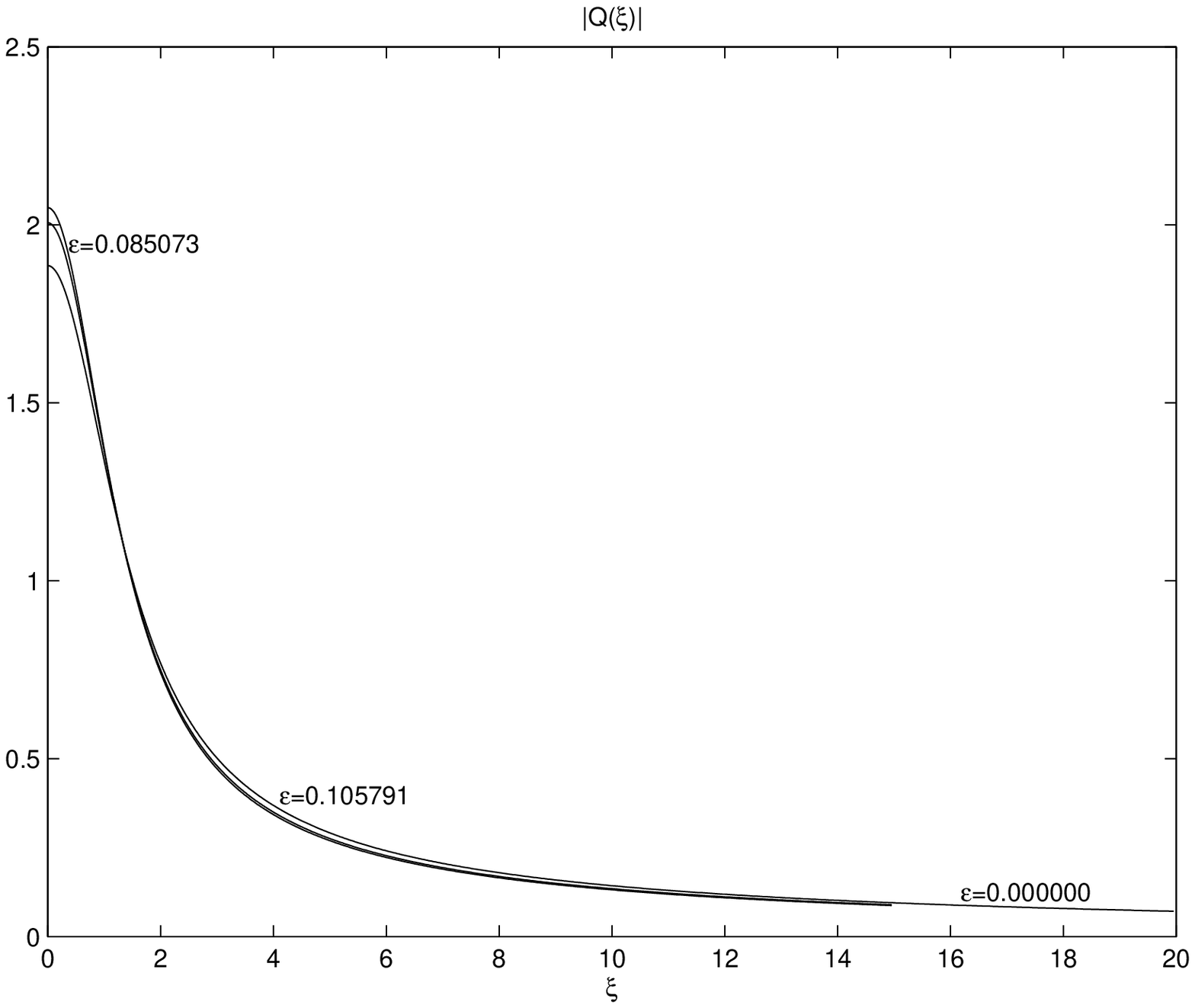, width=.98\linewidth}
\end{minipage}
\begin{minipage}{0.48\linewidth}
\centering\epsfig{file=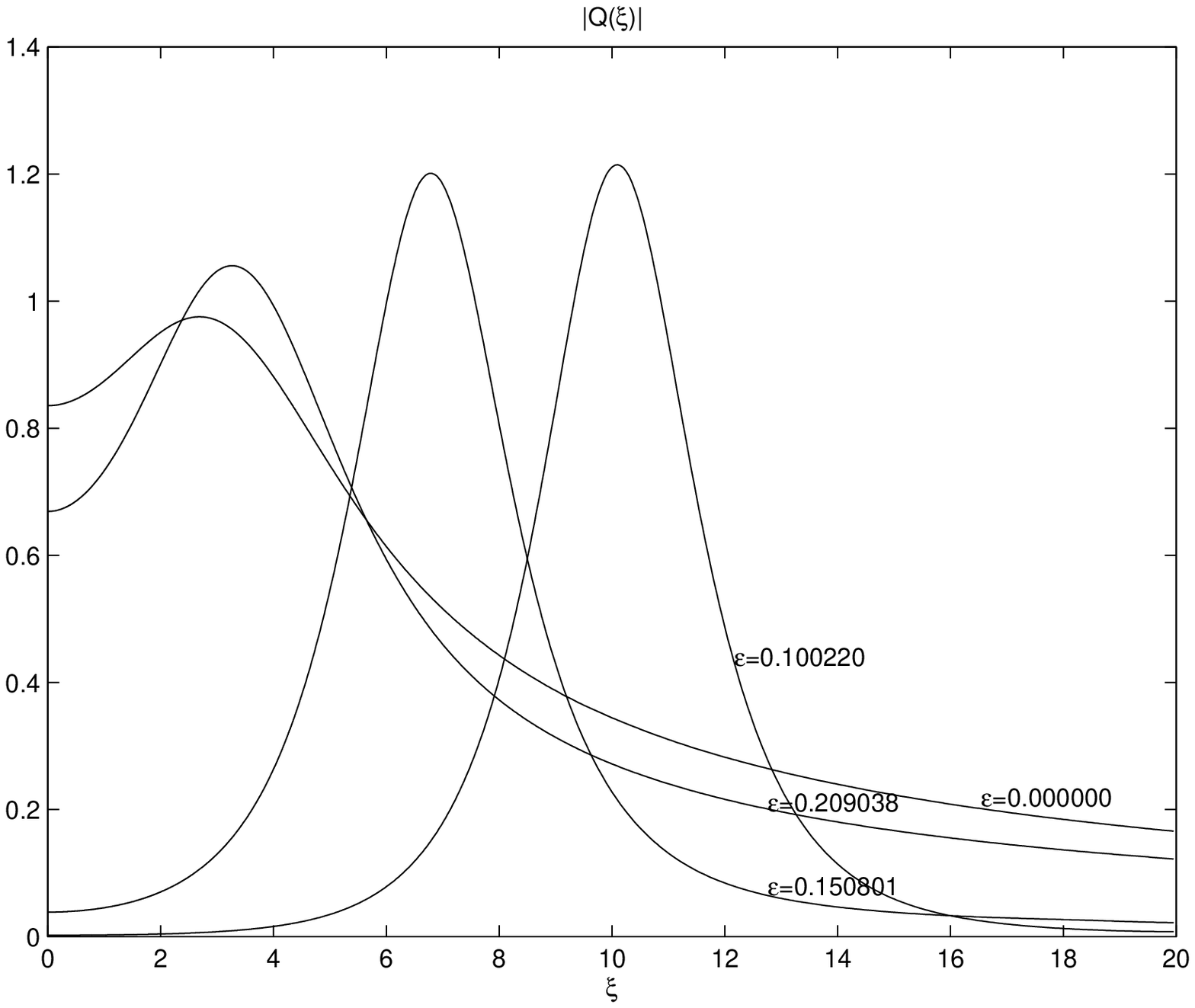, width=.98\linewidth}
\end{minipage}
\end{center}
\caption{{\sc Case II} ($d=3$, $\sigma=1$, $\delta=0$): the profile $|Q(\xi)|$ at points along
         the $j$-th branch. $j=1$ (left), $j=2$ (right); normalization
      $\omega=1$.}\label{figd3-b12}
\end{figure}

\begin{figure}[ht]
\begin{center}
\begin{minipage}{0.48\linewidth}
\centering\epsfig{file=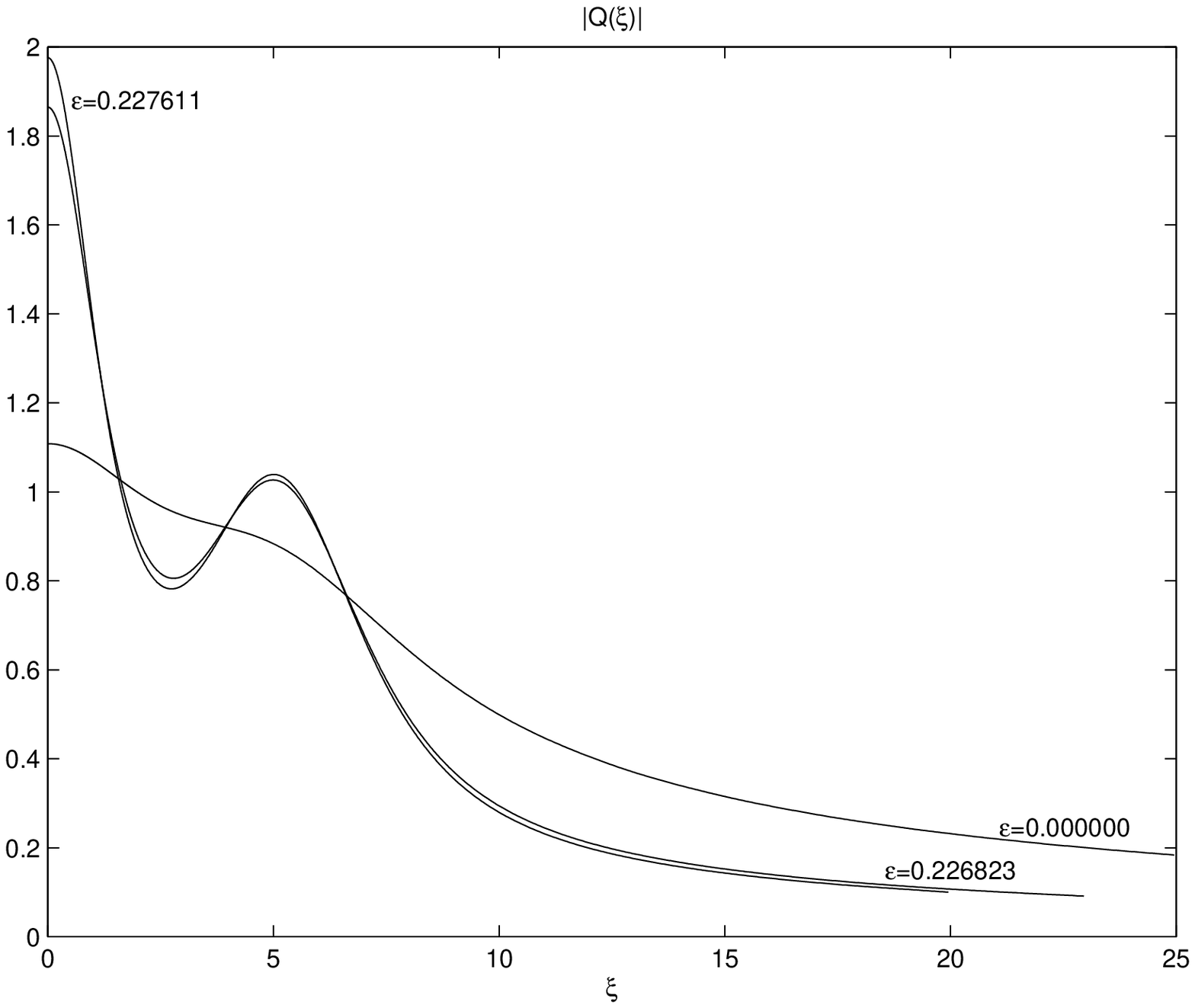, width=.98\linewidth}
\end{minipage}
\begin{minipage}{0.48\linewidth}
\centering\epsfig{file=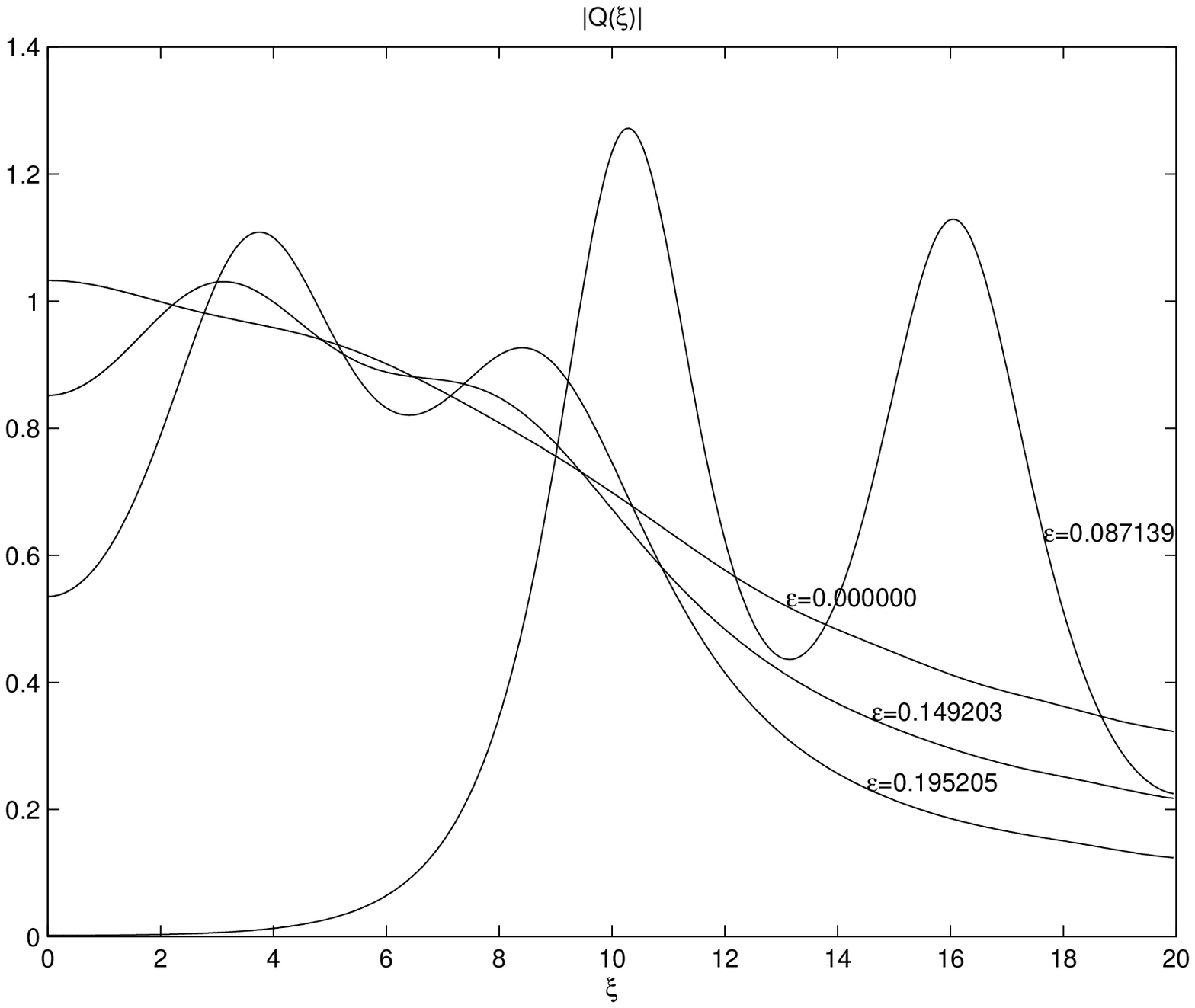, width=.98\linewidth}
\end{minipage}
\end{center}
\caption{{\sc Case II} ($d=3$, $\sigma=1$, $\delta=0$): the profile $|Q(\xi)|$ at points along
         the $j$-th branch $j=3$ (left), $j=4$ (right); normalization $\omega=1$.}\label{figd3-b34}
\end{figure}

\begin{remark}
From the regularity theory for CGL (see \cite{BartuccelliConstantinDoering})
 one expects
that the value of $\epsilon$ cannot cross a certain finite value
$\tilde\epsilon$, as we follow solutions along any branch. What we see
in the calculations 
is a turning point on each branch. In Table~\ref{tab1} and
Table~\ref{tab2}
 we list for each of the
calculated branches the parameter $\epsilon^*$ which is the maximal 
value of the perturbation parameter $\epsilon$  reached along 
the branch.
\end{remark}

\section{Stability of singular solutions}\label{sec4}

A natural question regarding the stability of the self-similar singularities constructed
in the previous section is for example the following: What will be the
behavior of
solutions to \VIZ{clg1a} and \VIZ{clg1b} if $u_0$ is a slight perturbation of a solution 
$Q$ of the profile equation \BVPQ\ and $f$ is small? (One can think for example of 
taking $u_0=\varphi Q+u_1$, where $\varphi$ is a compactly supported smooth cut-off function
which is identically one in a large ball,  $Q$ is a non-trivial solution of
\BVPQ, $u_1$ is small and compactly supported, and at the same time taking $f$ which
is small and compactly supported in $x$.)

Here we present an approach to this problem which  uses the method of dynamical rescaling
(see for example \cite{McLaughlinPapanicolaou}), which seems to be natural in this context.
We emphasize that we will not obtain fully rigorous analytical results which would
completely answer the question raised above. Our goal is to present some preliminary
calculations which seem to be adequate to the issue of stability in the
context of numerical simulations.
 
We briefly recall the main idea of the method of dynamical rescaling.
 We consider a solution $u:\RR{d}\times
(0,T)\to\CNUM$
 of the complex Ginzburg-Landau equation
\begin{equation}\label{clgS1}
\Im\PD{u}{t} + (1-\Im\epsilon) \LAP u + (1+\Im\delta) |u|^{2\sigma} u 
= 0\PERIOD
\end{equation}
Motivated by the scaling invariance of solutions
$$
u(x,t) \to \lambda^{1/\sigma} u(\lambda x,\lambda^2 t)
$$
of \VIZ{clgS1}, we write the solution $u$ as 
$$
u(x,t) = L^{-1/\sigma}(t) \, v\left(L^{-1}(t) x,\tau\right)\COMMA
$$
where $L(t)>0$ is to be chosen and $\d\tau = L^{-2}(t)\d
t$. From \VIZ{clgS1} we obtain
\begin{equation}\label{clgS2}
\Im\PD{v}{\tau} + \Im\kappa(\tau)\left(\xi\PD{v}{\xi} +
  \frac{1}{\sigma} v\right) + (1-\Im\epsilon)\LAP v +
(1+\Im\delta)|v|^{2\sigma} v = 0\COMMA
\end{equation}
where $\kappa(\tau) = - L^{-1}(\tau) \d L(\tau)/\d\tau = -L(t) \d
L(t)/\d t$.

We now chose $L(t)$ so that, roughly speaking, the typical length-scale over which $v$
oscillates is $1$. If $u$ develops a singularity  at $t=T$ (and it is
regular in $(0,T_1)$ for each $T_1<T$), then $L(t)>0$ in $(0,T)$
and $L(t)\to 0$ as $t\to T$.

One way to control oscillations of $v$ is to impose a condition
$\JFUN(v(.,\tau))=1$, where $\JFUN$ is a suitable functional
controlling regularity of solutions to \VIZ{clgS2}. One is then led to 
the following system
\begin{eqnarray} \label{clgS3}
&&\Im\PD{v}{\tau} + \Im\kappa(\tau)\left(\xi\PD{v}{\xi} +
  \frac{1}{\sigma} v\right) + (1-\Im\epsilon)\LAP v +
  (1+\Im\delta)|v|^{2\sigma} v = 0\COMMA \\
&& \JFUN(v(.,\tau)) = 1\COMMA
\end{eqnarray}
where the unknowns are $v$ and $\kappa(\tau)$.
Various choices
of the functional $\JFUN$ have been used for numerical
calculations. For example, in
\cite{McLaughlinPapanicolaou,LeMesurier} the functional 
$$
\JFUN_0(v) = \int_{\RR{d}} |\GRAD v|^2\d x
$$
was used. Although it does not directly follow from the known regularity
theory that $\JFUN_0$ controls the regularity of solutions to
\VIZ{clgS2} in the case $\epsilon=\delta=0$ this choice turned out to 
work satisfactorily in the numerical computations.

In the following analysis we assume that the functional $\JFUN$ is
invariant under the action of the symmetry group $S^1$, i.e.,
$\JFUN(\EXP{\Im\theta} v) = \JFUN(v)$ for all $\theta\in \PER$. 
Self-similar singularities of the form \VIZ{selfsim}
correspond to solutions of \VIZ{clgS3} of the form 
$v(\xi,\tau) =\EXP{\Im\omega\tau}Q(\xi)$, i.e., to $S^1$-orbits of solutions to the 
problem \BVPQ. We will consider the linearized stability of these orbits, which seems
to be the simplest natural notion of stability in our context.
Suppose $(v_0(\xi,\tau),\kappa_0(\tau))$ is a solution of \VIZ{clgS3} such 
that $v_0(\xi,\tau) = \EXP{\Im\omega_0\tau}Q(\xi)$ and
$\kappa_0(\tau)=\kappa_0=\CONST$. We
consider a perturbed solution $(v(\xi,\tau),\kappa(\tau))$ in the form
\begin{eqnarray*}
&& v(\xi,\tau) = \EXP{\Im\omega_0\tau}(Q(\xi) + w(\xi,\tau)) \\
&& \kappa(\tau) = \kappa_0 + \varkappa(\tau)\COMMA
\end{eqnarray*}
where $w$ and $\varkappa$ are infinitesimally small. A simple
calculation gives 
\begin{eqnarray}\label{clgS4}
&&\PD{w}{\tau} = \LOPER w + \varkappa(\tau) \left(\xi Q' +
  \frac{1}{\sigma}Q\right) \COMMA \\
&&\langle\DFUN{\JFUN}(Q),w\rangle = 0 \COMMA
\end{eqnarray}
where the operator $\LOPER$ is defined as the linearization of the
equation \VIZ{clgS2} along the solution $(v_0,\kappa_0)$, i.e.,
$$
\LOPER w = (\Im + \epsilon)\LAP w - \kappa_0
           \left(\xi\PD{w}{\xi} +\frac{1}{\sigma}w\right)
           -\Im\omega_0 w + \Im(1+\sigma) |Q|^{2\sigma} w + 
           \Im\sigma |Q|^{2\sigma-2}Q^2 \bar w
$$
and $\DFUN{\JFUN}$ denotes the  derivative of $\JFUN$. (Note
that the operator $\LOPER$ is not complex linear and therefore 
it is natural to carry out the analysis in the real representation.)

There are two eigenfunctions of $\LOPER$ that can be formally derived directly
from the invariance of the profile equation \VIZ{clg3} under the scaling symmetries
\begin{eqnarray}
(Q(\xi),\kappa,\omega) & \to & (\EXP{\Im\theta}Q(\xi),\kappa,\omega)
                            \COMMA\label{clgS5} \\
(Q(\xi),\kappa,\omega) & \to & (\lambda^{1/\sigma+\Im\omega/\kappa}Q(\lambda\xi),
                               \lambda^2\kappa,\lambda^2\omega) \PERIOD\label{clgS6} 
\end{eqnarray}\

The invariance under the $S^1$-symmetry \VIZ{clgS5} implies existence of an eigenfunction 
$Y_1 = \Im Q$ with the eigenvalue zero, while the symmetry under
\VIZ{clgS6} leads (formally) to an eigenfunction $Y_2 = \xi
Q' + (1/\sigma + \Im\omega_0/\kappa_0)Q$ with the eigenvalue
$2\kappa_0$. We note that $Y_2(\xi)=\BIGO{(\xi^{-1/\sigma-2})}$ as $\xi\to \infty$,
due to \VIZ{asymp1} and \VIZ{asymp2}. We have
\begin{eqnarray}
&& \LOPER Y_1 = 0 \COMMA \label{clgS7} \\
&& \LOPER Y_2 = 2\kappa_0 Y_2\PERIOD \label{clgS8}
\end{eqnarray} 

We rewrite \VIZ{clgS4} as 
\begin{eqnarray}\label{clgS9}
&& \PD{w}{\tau} = \LOPER w + \varkappa(\tau)\left(Y_2 -
               \frac{\omega_0}{\kappa_0} Y_1\right) \\
&& \langle\DFUN{\JFUN}(Q),w\rangle = 0\COMMA
\end{eqnarray}
where the unknowns are $w(\xi,\tau)$ and $\varkappa(\tau)$.

Proceeding further with our formal reasoning, we  view \VIZ{clgS9} as a 
dynamical system with one linear constraint in a suitable linear space $\WSP$ 
of functions on $\RR{d}$. In this formal analysis we will not try to
specify $\WSP$. For a rigorous analysis it would be natural  to try to find a suitable 
Banach space  containing  all smooth
compactly supported functions, together with the functions $Y_1$ and $Y_2$. 
Suppose that the space $\WSP$ can be decomposed
into a direct sum $\WSP = \RNUM Y_1 \oplus \RNUM Y_2 \oplus \ZSP$ where
$\ZSP\subset \WSP$ is invariant under
$\LOPER$. In the finite-dimensional situation, a sufficient condition for this would be that the
eigenvalues $0$ and $2\kappa_0$ are simple. Assuming that such a
decomposition exists, we write
\begin{equation}\label{clgS10}
w(\xi,\tau) = w_1(\tau) Y_1(\xi) + w_2(\tau) Y_2(\xi) + \ZFUN(\xi,\tau)\COMMA
\end{equation}
where $w_1$, $w_2$ are scalar functions and $\ZFUN(.,\tau)\in\ZSP$. 

After substituting \VIZ{clgS10}
into \VIZ{clgS9} and using  $\langle \DFUN{\JFUN}(Q),Y_1\rangle = 0$ (a consequence
of the invariance of $\JFUN$), we obtain
\begin{eqnarray}\label{clgS11}
\PD{\ZFUN}{\tau} & = & \LOPER \ZFUN \\
\PD{w_1}{\tau} & = & -\varkappa(\tau) \frac{\omega_0}{\kappa_0}  \\
\PD{w_2}{\tau} & = & 2\kappa_0 w_2 + \varkappa(\tau) \\
w_2 & = & -\frac{\langle\DFUN{\JFUN}(Q),\ZFUN\rangle}%
                           {\langle\DFUN{\JFUN}(Q),Y_2\rangle}\label{clgS11a}
\end{eqnarray}
where we assume that $\langle\DFUN{\JFUN}(Q),Y_2\rangle\neq 0$. This
 is a completely natural assumption in the context of
\VIZ{clgS3}. From \VIZ{clgS11} - \VIZ{clgS11a} one can see that in a finite-dimensional
situation and under the assumptions 
stated above, the following conditions would be equivalent:
\begin{description}
\item[{\rm (i)}] If $(w,\varkappa)$ is a solution of \VIZ{clgS4}, then 
  $w(\xi,\tau)$ approaches exponentially $a Y_1(\xi)$ for some $a\in\RNUM$ and
  $\varkappa(\tau)$ approaches exponentially zero.
\item[{\rm (ii)}] The spectrum of  $\LOPER|_{\ZSP}$ belongs to the set 
$\{z \in\CNUM\SEP \REAL{z} < 0\}$.
\end{description}

In a finite-dimensional situation and under the assumption that the eigenvalues 
$0$ and $2\kappa_0$ are simple,
condition (ii) would be  equivalent to the
condition that all the spectrum of $\LOPER$ except $0$ and $2\kappa_0$
lies in $\{z \in\CNUM\SEP \REAL{z} < 0\}$.

The condition (i) is exactly the linearized orbital stability of the
orbit $\left(\EXP{\Im\omega_0\tau}Q,\kappa_0\right)$ for the
system \VIZ{clgS3}, \VIZ{clgS4}.

It seems to be a non-trivial problem to put the above formal analysis on a
rigorous basis in the infinite-dimensional setting. However, the formal
analysis strongly suggests that 
in a finite dimensional situation which arises in numerical approximations
of \VIZ{clgS3} the stability of the solutions 
$\left(\EXP{\Im\omega_0\tau}Q(\xi),\kappa_0\right)$
should be governed by the spectrum of an appropriate approximation of
the operator $\LOPER$. 

We remark that our stability analysis is
independent of $\JFUN$, except for the natural assumption
$\langle\DFUN{\JFUN}(Q),Y_2\rangle\neq 0$.

\medskip

We numerically calculated the approximation of eigenvalues for a
discrete approximation of $\LOPER$
{\it in the space of radial functions} along the branches
of solutions parameterized by $\epsilon$. 
The original system \VIZ{clgS3} was truncated
to a finite interval $(0,\xi_1)$ by imposing the time-dependent boundary condition 
at $\xi=\xi_1$
\begin{equation}\label{clgS12}
\PD{v}{\tau}(\xi_1,\tau) + \kappa(\tau)\xi_1\PD{v}{\xi}(\xi_1,\tau) +
\frac{\kappa(\tau)}{\sigma} v(\xi_1,\tau) = 0\COMMA
\end{equation}
which is a time-dependent equivalent of the boundary
condition derived from the asymptotic expansion in
Section~\ref{sec2}. The condition \VIZ{clgS12} was also used in the
numerical simulations in \cite{Zacharov2} and \cite{LandmanPapanicolaou}.
The computations were carried out 
for $\xi_1=30$.  

The accuracy of our numerical approximation can be checked indirectly
by comparing the predicted eigenvalues $0$ and $2\kappa_0$ with the
corresponding eigenvalues we obtained from our numerical calculations.
We saw a very good agreement, in most cases the error was of the order
$10^{-4}$.
 Both $0$ and $2\kappa_0$ appeared simple, except in some natural 
degenerate cases
when other eigenvalues were crossing them as we moved
along branches. 

\noindent
{\bf Case I} ($d=1$, $\sigma=2.3$, $\delta=0$): The calculations
confirmed what one intuitively expects:
The solutions on the upper part of the branch $j=1$ are 
stable and all other solutions are unstable. In Figure~\ref{fig1} we
used a solid line for the stable parts of the curves and a dashed line
for the unstable parts. \\

\noindent
{\bf Case II} ($d=3$, $\sigma=1$, $\delta=0$): The situation is
similar with one notable exception. In our computations we detected
stable
solutions also on the lower part of the branch $j=2$. Accuracy of our
approximation did not allow us to decide whether all solutions on the
lower part of the branch $j=2$ for $0<\epsilon<\epsilon^*$ are
stable, since as $\epsilon$ approached zero we observed some 
eigenvalues very close to the imaginary axis. In Figure~\ref{fig2} we
used a solid line for the stable parts of the curves and a dashed line
for the unstable parts.\\

We  emphasize again that the calculations we
carried out only deal with stability in the space of radial functions. 
Based on computations in \cite{LandmanPapanicolaou}, it appears that the
solutions on the upper part of the branch $j=1$ are also stable with
respect to perturbations that break the radial symmetry. 

It is clear that further work is required to fully clarify issues
concerning the stability of  solutions described in this paper.

\section{Appendix}
In this section we give a full proof of Theorem~\ref{thm1}.

We fix  $r_1>0$ and $\vartheta>0$.
Values of these parameters will be chosen later. We denote by
$\BSPT$ the Banach space of continuous functions
$u:[\xi_1,\infty)\to\CNUM$ for which the norm
$$
\NORM{u}{\vartheta} = \sup_{\xi\geq\xi_1} |\xi|^{1/\sigma-\vartheta}
|u(\xi)|
$$
is finite. In this proof we use the following notation
\begin{eqnarray*}
&& \LAMT =  [\kappa_1,\kappa_2]\times [\omega_1,\omega_2] \times
[0,\epsilon_1]\times [\delta_1,\delta_2] \\
&& \LAM = \{\gamma\in\CNUM\SEP |\gamma|\leq r_1\}\times \LAMT \PERIOD
\end{eqnarray*}
A point in $\LAMT$ is denoted by
$\lambdat=(\lambda_1,\lambda_2,\lambda_3,\lambda_4) \equiv
(\kappa,\omega,\epsilon,\delta)$ and similarly $\lambda\in\LAM$
represents
$(\lambda_0,\lambda_1,\lambda_2,\lambda_3,\lambda_4) \equiv
(\gamma,\kappa,\omega,\epsilon,\delta)$. The functions $P$, $E$, $W$,
and $K$ introduced in Section~\ref{sec2} are written as
$P=P(\lambdat,\xi)$, 
$E=E(\lambdat,\xi)$, 
$W=W(\lambdat,\xi,\eta)$, 
$K=K(\lambdat,\xi)$. 
With a slight abuse of notation we sometimes write
also $P=P(\lambda,\xi)$ even if the function does not depend on
$\lambda_0$.

\begin{lemma} \label{lem1}
There exists $C>0$ such that for $\lambdat\in\LAMT$,
$\xi,\eta\geq\xi_0$ we have
\begin{equation}\label{clgA1}
|K(\lambdat,\xi,\eta)| \leq \left\{\begin{array}{ll}
              {\displaystyle C\xi^{-1/\sigma}\eta^{1/\sigma-1}} &
              \mbox{ for }\; \xi_0\leq\eta\leq \xi\COMMA \\
              {\displaystyle C\xi^{-d+1/\sigma}\eta^{-1-1/\sigma+d}} &
              \mbox{ for }\; \xi_0\leq\xi\leq \eta\PERIOD
              \end{array}\right.
\end{equation}
\end{lemma}

\PROOF This statement is an easy consequence of the definition
of $K$ and of the asymptotic expansion for $U$ as $z\to\infty$ in the
sector $-\pi/2\leq\ARG{z}\leq\pi$. An important point in connection with 
\VIZ{clg9} is that the constant in the remainder $\BIGO{(|z|^{-n})}$
can be taken same when the parameter $a$ runs through a compact subset 
of $\{ z\SEP \REAL{z} > 0\}$. (This can be seen for example from \VIZ{clg8}.)\qed

\begin{lemma}\label{lem2}
  Assume $(2\sigma+1)\vartheta < 2 + 2/\sigma -d$. Then the formula
  \begin{equation}\label{clgA2}
     T(\lambda,u)(\xi) = \gamma P(\lambdat,\xi) - \int_{\xi_1}^\infty
     (1+\Im\delta) K(\lambdat,\xi,\eta) |u(\eta)|^{2\sigma}u(\eta)\d\eta
   \end{equation}
  defines a continuous mapping $T:\LAM\times \BSPT \to
  \BSPT$. 

  Moreover, there exists $C>0$ such that
  \begin{equation}\label{clgA3}
    \NORM{T(\lambda,u)}{\vartheta} \leq C |\gamma|\xi_1^{-\vartheta} + C
    \xi_1^{-2+2\sigma\vartheta} \NORM{u}{\vartheta}^{2\sigma+1}
  \end{equation}
  and
  \begin{equation}\label{clgA4}
    \NORM{T(\lambda,u)-T(\lambda,v)}{\vartheta} \leq C
    \xi_1^{-2+2\sigma\vartheta} 
    \NORM{u-v}{\vartheta} \left(\NORM{u}{\vartheta}^{2\sigma} +
    \NORM{v}{\vartheta}^{2\sigma}\right) 
  \end{equation}
  for all $\lambda\in\LAM$ and all $u,v\in\BSPT$.
\end{lemma}

\PROOF The convergence of the integrals in \VIZ{clgA2} and
\VIZ{clgA3} follows easily from Lemma~\ref{lem1}. To get the estimate
\VIZ{clgA4} we use Lemma~\ref{lem1} together with the elementary
inequality 
$$
\left| |z_1|^{2\sigma}z_1 -  |z_2|^{2\sigma}z_2\right| \leq M
|z_1-z_2| \left(|z_1|^{2\sigma} + |z_2|^{2\sigma}\right) \COMMA
$$
which holds for all $z_1, z_2\in\CNUM$ if $M>0$ is properly chosen.

The continuity of the mapping $T$ can be proved by suitably splitting the
integral over $(\xi_1,\infty)$ into two integrals, one over $(\xi_1,\xi_2)$
and the other one over
$(\xi_2,\infty)$. We can then estimate the integral over $(\xi_2,\infty)$ by
using Lemma~\ref{lem1}. \qed

The estimates \VIZ{clgA3} and \VIZ{clgA4} show that $T(\lambda,.)$ is
a contraction of the ball $\BSPB_\rho=\{u\in\BSPT\SEP
\NORM{u}{\vartheta} \leq\rho\}$ into itself if $\rho$ is such that
\begin{eqnarray}
&& Cr_1\xi_1^{-\vartheta} + C
\xi_1^{-2+2\sigma\vartheta}\rho^{2\sigma+1} \leq
\rho\COMMA\label{clgA5} \\
&& 2C\rho^{2\sigma}\xi_1^{-2+2\sigma\vartheta} < 1\PERIOD \label{clgA6}
\end{eqnarray}

\begin{proposition}\label{prop1}
Assume $(2\sigma+1)\vartheta < 2 + 2/\sigma -d$. If \VIZ{clgA5} and
\VIZ{clgA6} are satisfied, then for each $\lambda\in\LAM$ the mapping
$T(\lambda,.)$ has a unique fixed point $u(\lambda)=u(\lambda,\xi)$ in
$\BSPB_\rho=\{u\in\BSPT\SEP \NORM{u}{\vartheta}
\leq\rho\}$.

Moreover, the mapping $\lambda\mapsto u(\lambda)$ is a continuous
mapping from $\LAM$ to $\BSPT$.
\end{proposition}

\PROOF The proposition follows directly from Lemma~\ref{lem2}
and the Banach Fixed Point Theorem. \qed

Next we need to establish that the solutions $u(\lambda,\xi)$ have the 
required regularity properties at the infinity.

\begin{lemma}\label{lem3}
The function $u(\lambda)=u(\lambda,\xi)$ from Proposition~\ref{prop1}
exhibits the following behavior as $\xi\to\infty$:
\begin{equation}\label{clgA7}
u(\lambda,\xi) = \xi^{-1/\sigma - \Im\omega/\kappa}
\left(\sum_{j=0}^n a_j(\lambda) \xi^{-2j} +
  \BIGO{(\xi^{-2n-2})}\right)\COMMA
\end{equation}
where $a_j(\lambda)$ are continuous functions of $\lambda$ and the
constant in $\BIGO{(\xi^{-2n-2})}$ is independent of
$\lambda\in\LAM$. 

Moreover, the expansion \VIZ{clgA7} can be differentiated in the
following sense:
\begin{equation}\label{clgA8}
\frac{\partial^k u}{\partial\xi^k}(\lambda,\xi) = 
\frac{\partial^k}{\partial\xi^k}\left(\xi^{-1/\sigma-\Im\omega/\kappa} 
  \sum_{j=0}^n a_j(\lambda) \xi^{-2j} \right) +
\BIGO{(\xi^{-1/\sigma-2n-2-k})}\COMMA
\end{equation}
where the constant in $\BIGO{(\xi^{-1/\sigma-2n-2-k})}$ is independent of
$\lambda\in\LAM$.
\end{lemma}

\PROOF The proof is similar to bootsrapping arguments used
in the regularity theory. From \VIZ{clgA2} and the assumption
$(2\sigma+1)\vartheta< 2 + 2/\sigma -d$
we get immediately
\begin{equation}\label{clgA9}
u(\lambda,\xi) = \BIGO{(\xi^{-1/\sigma})}\;\;\mbox{as}\;\;
\xi\to\infty\PERIOD
\end{equation}
We now rewrite the equation 
$u= T u$ 
as
\begin{equation}\label{clgAR}
u(\xi)=(\gamma+\gamma_1)P(\xi)
-\frac{1+\Im\delta}{1-\Im\epsilon}  P(\xi)
\int_{\xi}^{\infty} EW^{-1}|u|^{2\sigma}u+
\frac{1+\Im\delta}{1-\Im\epsilon}E(\xi)
\int_{\xi}^\infty PW^{-1}|u|^{2\sigma}u\COMMA
\end{equation}
where 
$$
\gamma_1=\frac{1+\Im\delta}{1-\Im\epsilon}
\int_{\xi_1}^{\infty} EW^{-1}|u|^{2\sigma}u\PERIOD
$$
Using \VIZ{clgA9} and \VIZ{clgAR}, we  obtain (for a suitable $a_0$)
\begin{equation}\label{A2}
u(\xi)=\xi^{-1/\sigma-\Im\omega/\kappa}\left(a_0+\BIGO(\xi^{-2})\right)
\end{equation}
We now repeat the procedure and use \VIZ{clgAR} with \VIZ{A2} instead of \VIZ{clgA9}.
After integrating by parts in integrals of the form
$$
\int_\xi^{\infty}\eta^\alpha \EXP{-\nu\eta^2/2}=
-\frac1\nu\int_\xi^{\infty}\eta^{\alpha-1}\frac{\partial}{\partial\eta}\EXP{-\nu\eta^2/2}
$$
which come up in the calculation, we get easily
$$
u(\xi)=\xi^{-1/\sigma-\Im\omega/\kappa}
\left(a_0+a_1\xi^{-2}+\BIGO(\xi^{-4})\right)\PERIOD
$$
Repeating this procedure, we get \VIZ{clgA7} by induction.
The expansion \VIZ{clgA8} can be obtained in a similar way, after we differentiate
\VIZ{clgAR}.
\qed

\begin{lemma}\label{lem4}
The mapping $\lambda\mapsto u(\lambda)$ from Proposition~\ref{prop1}
is continuously differentiable as a mapping from $\LAM$ to
$\BSPT$ 
\end{lemma}

\PROOF We define $\phi:\LAM\times\BSPT
\to\BSPT$ as $\phi(\lambda,v) = v - T(\lambda,v)$ and we
denote $\Mm=\{ (\lambda,u(\lambda)),\lambda\in\LAM\}\subset
\LAM\times\BSPT$. By using similar estimates as those used
in the proof of Lemma~\ref{lem2}, one can easily see that $\phi$ is
differentiable with respect to $v$ and that the (partial) derivative
$\PDER{v}{\phi}$ is continuous as a map from
$\LAM\times\BSPT$ into the space $\Ll(\BSPT,\BSPT)$ of bounded linear
operators on $\BSPT$. Moreover, a calculation similar to the one leading
to \VIZ{clgA4} gives
\begin{equation}\label{A3}
\NORM{I-\PDER{v}{\phi}}{}\le 2C\xi_1^{-2+2\sigma\vartheta}\rho^{2\sigma}
\end{equation}
in $\LAM\times\BSPB_\rho$.
Therefore, since we assume \VIZ{clgA6}, 
$\PDER{v}{\phi}$ is invertible
for all $(\lambda,v)\in\LAM\times\BSPB_\rho$. By Lemma~\ref{lem5}
below, $\phi$ is also differentiable with respect to $\lambda$ at each 
point of $\Mm$ and the partial derivative
$\PDER{\lambda}{\phi}(\lambda,v)$ is continuous in $\Mm$. Now we can
conclude with the same arguments as in the Implicit Function Theorem
to show that $u$ is differentiable and
\begin{equation}\label{clgAder}
\PD{u}{\lambda}(\lambda) = -
\left[\PDER{v}{\phi}(\lambda,u(\lambda))\right]^{-1}\cdot
\left[\PDER{\lambda}{\phi}(\lambda,u(\lambda))\right]\PERIOD 
\end{equation}
\qed

Lemma~\ref{lem4} also implies the differentiability of the function
$\lambda\mapsto \PD{u}{\xi}(\lambda,\xi_1)$, since the partial
derivative of $u$ at $\xi_1$ can be expressed in terms of $u$ due to
the fact that $u$ solves a differential equation.
For example, one can take a smooth function $\varphi:[\xi_1,\infty)
\to\RNUM$  with $\varphi(\xi_1)=1$ and $\varphi(\xi)\equiv 0$, for
$\xi>\xi_1+1$ and write $u(\xi_1) = -\int_{\xi_1}^\infty
(u\varphi)''\d\xi$. 
Then we express $u''$ from the differential equation
for the profile \VIZ{clg3}
and integrate by parts to obtain
$$
\PD{u}{\xi}(\xi_1) = c_1 u(\xi_1) + \int_{\xi_1}^\infty
(u(\xi)\psi(\xi) + |u(\xi)|^{2\sigma}u(\xi)\tilde\psi(\xi)) \d\xi\COMMA
$$
where $\psi$, $\tilde\psi$ are supported on a finite interval and $c_1$ is a
suitable constant.
\qed

\begin{lemma}\label{lem5}
The mapping $\phi:\LAM\times\BSPT
\to\BSPT$ defined as $\phi(\lambda,v) = v - T(\lambda,v)$ is
differentiable with respect to $\lambda$ at each point 
$\Mm=\{ (\lambda,u(\lambda)),\lambda\in\LAM\}\subset
\LAM\times\BSPT$. Moreover, the partial derivative
$\PDER{\lambda}{\phi}$ is continuous in $\Mm$.
\end{lemma}

\PROOF Using \VIZ{clg6} we note that
\begin{equation}\label{clgA16}
\PD{U}{a}(a,b,z) = z^{-a}\log z \left(1 + \BIGO(|z|^{-1})\right)\COMMA
\end{equation}
as $z\to\infty$ in the sector $-\pi/2\leq\ARG{z}<\pi$ and $a$
belongs to a compact subset of
$\{\zeta\SEP\REAL{\zeta}>0\}\subset\CNUM$. 
Let us carry out the proof for the partial derivative $\partial\phi/\partial\lambda_1$,
for example. The other partial derivatives can be handled in a similar way. (If fact,
for $\lambda_0$,$\lambda_2$ and $\lambda_4$ the proof is much easier.)
When $\epsilon>0$ one sees easily (using the rapid decay of $W^{-1}(\eta)$ in that
case) that $\phi$ is differentiable in $\lambda_1$ and that 
\begin{equation}\label{A4}
\PD{\phi}{\lambda_1}(\lambda,v)(\xi) = \gamma\PD{P}{\lambda_1}(\lambda,\xi) 
+ \int_{\xi_1}^\infty \PD{K}{\lambda_1}(\lambda,\xi,\eta) \,
|v(\eta)|^{2\sigma} v(\eta)\d\eta\PERIOD
\end{equation}

It follows from \VIZ{clgA16} that, for any $\vartheta>0$,
 the first term on the right-hand side is 
continuous as a mapping from $\LAM$ to $\BSPT$. 
For any $\epsilon_2>0$, the integral on the right-hand side has the required 
continuity properties in the region 
$\epsilon_2\le\epsilon\le \epsilon_1$.
due to the rapid decay of $W^{-1}(\eta)$. The problem is to obtain estimates
which are uniform for $0<\epsilon<\epsilon_1$.
A closer inspection of
the integral on the right-hand side of \VIZ{A4} reveals that the only terms for
which the required continuity is not obvious come from differentiating 
the exponential
$\EXP{-\Im\kappa\eta^2/2(1-\Im\epsilon)}$ that appears in
 $W^{-1}(\eta)$
This term can be handled with using
integration by parts as follows.
  Recalling that $\lambda_1=\kappa$, we
must estimate an integral of the form
\begin{equation}\label{clgA17}
\int_\xi^\infty P(\lambda,\eta) \eta^{d-1}
\frac{\partial}{\partial\kappa} \left(
  \EXP{\frac{\Im\kappa}{1-\Im\epsilon}\frac{\eta^2}{2}}\right)
|v(\eta)|^{2\sigma} v(\eta) \d\eta \PERIOD
\end{equation}
To estimate this integral we write
$$
\frac{\partial}{\partial\kappa}\left(
  \EXP{\frac{\Im\kappa}{1-\Im\epsilon}\frac{\eta^2}{2}}\right) = 
\frac{\eta}{2\kappa}\frac{\partial}{\partial\eta}
\left(\EXP{\frac{\Im\kappa}{1-\Im\epsilon}\frac{\eta^2}{2}}\right)
$$
and integrate by parts. This  eliminates the power
$\eta^2$ obtained after differentiating in
\VIZ{clgA17}, and the required estimate follows.
 \qed

\medskip

To complete the proof of Theorem~\ref{thm1} we need to show that the coordinate
$\gamma=\lambda_0$ can be replaced by $\beta=u(\lambda,\xi_1)$, if $r_1$ and
$\rho$ are chosen properly. We have
\begin{equation}\label{A5}
\frac{\partial u}{\partial\lambda_0}(\lambda,\xi_1)=
P(\lambda,\xi_1)+\int_{\xi_1}^\infty K(\tilde\lambda,\xi_1,\eta)
\frac{\partial }{\partial\lambda_0} |u|^{2\sigma}u(\lambda,\eta)\d\eta 
\PERIOD
\end{equation}
(Since $\lambda_0$ is complex, we 
interpret this, with a slight abuse of notation, as an equation between real 
$2\times 2$ matrices. Another possibility would be to interpret it literally and
do a similar calculation for ${\partial}/{\partial\bar\lambda_0}$.)
We carry out the differentiation in the integral on the right-hand side
and use \VIZ{clgAder}, which gives
\begin{equation}\label{A6}
\PD{u}{\lambda_0}(\lambda) = -
\left[\PDER{v}{\phi}(\lambda,u(\lambda))\right]^{-1}\cdot P\PERIOD 
\end{equation}
We note from \VIZ{A3} that by a suitable choice of the constant $\rho$ 
(to be specified below) we can achieve that 
\begin{equation}\label{A7}
\NORM{[\PDER{v}{\phi}(\lambda,u(\lambda))]^{-1}}{}\le 2\COMMA
\end{equation} 
in $\LAM$
where the norm is taken in the space of linear operators on the space
$\mathcal{X}_0$.
Using \VIZ{A5} - \VIZ{A7}, we see easily, by a similar calculation as in the
proof of \VIZ{clgA4}, that by a suitable choice of $\rho$ we can achieve
that 
$\frac{\partial u}{\partial\lambda_0}(\lambda,\xi_1)/P(\lambda,\xi_1)$ 
is close to the identity uniformly
in $\LAM$. A suitable choice for the constants $\rho$ and $r_1$
is, for example, as follows
$$
\rho=\alpha\xi_1^{1/\sigma-\vartheta}, \quad 
\quad r_1=\alpha^{2\sigma+1}\xi_1^{1/\sigma}\COMMA
$$
where $\alpha$ is sufficiently small.
The proof of Theorem~\ref{thm1} can now be easily completed.
\qed

\bibliographystyle{siam}
\bibliography{clg}

\end{document}